\documentclass[reqno, draft]{amsart}
\usepackage{amssymb}
\usepackage{upref}
\usepackage{graphicx}

\newcommand{\bbN}{{\mathbb{N}}}
\newcommand{\bbR}{{\mathbb{R}}}

\newcommand{\bbZ}{{\mathbb{Z}}}
\newcommand{\bbC}{{\mathbb{C}}}
\newcommand{\bbT}{{\mathbb{T}}}

\newcommand{\calB}{{\mathcal B}}

\newcommand{\calC}{{\mathcal C}}

\newcommand{\calD}{{\mathcal D}}

\newcommand{\calM}{{\mathcal M}}
\newcommand{\calE}{{\mathcal E}}
\newcommand{\calL}{{\mathcal L}}
\newcommand{\calF}{{\mathcal F}}
\newcommand{\calK}{{\mathcal K}}

\newcommand{\N}{g}

\newcommand{\dott}{\,\cdot\,}
\newcommand{\hatt}{\widehat}  

\newcommand{\Div}{\operatorname{Div}}

\newcommand{\no}{\nonumber}
\newcommand{\lb}{\label}
\newcommand{\f}{\frac}
\newcommand{\ul}{\underline}
\newcommand{\ol}{\overline}
\newcommand{\Oh}{O}

\newcommand{\humu}{{ \hat{\underline{\mu} }}}
\newcommand{\hunu}{{\underline{\hat{\nu}}}}
\newcommand{\hmu}{{\hat{\mu} }}
\newcommand{\hnu}{{\hat{\nu}}}

\newcommand{\uz}{{\underline{z}}}

\newcommand{\uxi}{{\underline{\Xi}}}

\newcommand{\ual}{{\underline{\alpha}}}
\newcommand{\ua}{{\underline{A}}}

\newcommand{\Pinfp}{{P_{\infty_+}}}
\newcommand{\Pinfm}{{P_{\infty_-}}}
\newcommand{\Pzp}{{P_{0,+}}}
\newcommand{\Pzm}{{P_{0,-}}}
\newcommand{\Pinfpm}{{P_{\infty_\pm}}}
\newcommand{\Pzpm}{{P_{0,\pm}}}
\newcommand{\Pinfmp}{{P_{\infty_\mp}}}

\renewcommand{\Re}{\text{\rm Re}}
\renewcommand{\Im}{\text{\rm Im}}

\DeclareMathOperator{\sym}{Sym}

\DeclareMathOperator{\sSB}{s-SB}

\allowdisplaybreaks
\numberwithin{equation}{section}

\newtheorem{theorem}{Theorem}[section]

\newtheorem{lemma}[theorem]{Lemma}
\newtheorem{corollary}[theorem]{Corollary}

\theoremstyle{definition}

\newtheorem{remark}[theorem]{Remark}

\newtheorem{example}[theorem]{Example}

\newcommand{\abs}[1]{\lvert#1\rvert}

\begin{document}
\title[Algebro-Geometric Solutions of a Discrete System]{Algebro-Geometric
Solutions of a Discrete System Related to the Trigonometric Moment
Problem}
\author[J.\ Geronimo]{Jeffrey S.\ Geronimo}
\address{School of Mathematics, Georgia Institute of Technology,
Atlanta, GA 30332-01660, USA}
\email{geronimo@math.gatech.edu}
\urladdr{http://www.math.gatech.edu/\~{}geronimo/}
\author[F.\ Gesztesy]{Fritz Gesztesy}
\address{Department of Mathematics,
University of Missouri,
Columbia, MO 65211, USA}
\email{fritz@math.missouri.edu}
\urladdr{http://www.math.missouri.edu/people/fgesztesy.html}
\author[H.\ Holden]{Helge Holden}
\address{Department of Mathematical Sciences,
Norwegian University of
Science and Technology, NO--7491 Trondheim, Norway}
\email{holden@math.ntnu.no}
\urladdr{http://www.math.ntnu.no/\~{}holden/}
\thanks{The research of the second and third author was supported in part by
the Research Council of Norway.}
\thanks{Supported in part by the US National Science
Foundation under Grants No.\ DMS-0200219 and DMS-0405526.}
\date{July 23, 2004}

\begin{abstract}
We derive theta function representations of algebro-geometric solutions 
of a discrete system governed by a transfer matrix associated with (an
extension of) the trigonometric moment problem studied by Szeg{\H o}
and Baxter. We also derive a new hierarchy of coupled nonlinear difference
equations satisfied by these algebro-geometric solutions.
\end{abstract}

\maketitle

\section{Introduction}\lb{sbs1}

Let $\{\alpha (n)\}_{n\in\bbN}\subset\bbC$ be a sequence of complex
numbers subject to the condition
\begin{equation}
       |\alpha (n)|<1 \, \text{  for all $n\in\bbN$,} \lb{sb1.0}
\end{equation}
and define the transfer matrix
\begin{equation}
T(z)=\begin{pmatrix} z & \alpha \\ \overline{\alpha} z &1
\end{pmatrix}, \quad z\in\bbT, \lb{sb1.3}
\end{equation}
with spectral parameter $z$ on the unit circle
$\bbT=\{z\in\bbC\,|\,|z|=1\}$. Consider the system of difference
equations
\begin{equation}
\Phi(z,n) = T(z,n)\Phi(z,n-1), \quad
(z,n)\in\bbT\times\bbN  \lb{sb1.1}
\end{equation}
with initial condition $\Phi (z,0) = \left(\begin{smallmatrix} 1\\ 1
\end{smallmatrix}\right)$, 
$z\in\bbT$, where
\begin{equation}
\Phi(z,n)=\begin{pmatrix} \varphi(z,n)\\ \varphi^*(z,n) \end{pmatrix},
\quad (z,n)\in\bbT\times\bbN_0.
\lb{sb1.4}
\end{equation}
(Here $\bbN_0=\bbN\cup\{0\}$.) Then $\varphi (\dott,n)$ are monic
polynomials of degree $n$ and
\begin{equation}
\varphi^* (z,n)= z^n{\ol \varphi}(1/z,n), \quad (z,n)\in\bbT\times\bbN_0,
\lb{sb1.5}
\end{equation}
the reversed ${}^*$-polynomial of $\varphi(z,n)$, is of degree at
most $n$. These polynomials were first introduced by Szeg\H o in the
1920's in his work on the asymptotic distribution of eigenvalues of
sections of Toeplitz forms \cite{Szego:1920}, \cite{Szego:1921}  (see
also
\cite[Chs.\ 1--4]{GrenanderSzego:1984},
\cite[Ch.\ XI]{Szego:1978}). Szeg\H o's point of departure was the
trigonometric moment problem and hence the theory of orthogonal
polynomials on the unit circle: Given a probability measure $d\sigma$ 
supported on an infinite set on the unit circle, find monic polynomials of
degree $n$ in $z=e^{i\theta}$, $\theta\in [0,2\pi]$, such that
\begin{equation}
\int_{0}^{2\pi} \gamma(n)^2 d\sigma(e^{i\theta}) \, 
\overline{\varphi (e^{i\theta},m)} \varphi (e^{i\theta},n)
=\delta_{m,n}, \quad m,n\in\bbN_0,
\lb{sb1.8}
\end{equation}
where
\begin{equation}
\gamma(n)^2=\begin{cases} 1& \text{for $n=0$,} \\
\prod_{j=1}^n \big(1-|\alpha(j)|^2\big)^{-1} & \text{for $n\in\bbN$.}
\end{cases} \lb{sb1.7}
\end{equation}
Here we chose to emphasize monic polynomials
$\varphi(\dott,n)$ in order to keep the factor $\gamma$ out of the transfer
matrix $T$. Szeg\H o showed that the polynomials \eqref{sb1.4} satisfy
the recurrence formula \eqref{sb1.1}.  Early work in this area includes
contributions by Akhiezer \cite[Ch.\ 5]{Akhiezer:1965},
Geronimus \cite{Geronimus:1946}, \cite{Geronimus:1948},
\cite[Ch.\ I]{Geronimus:1961}, Krein \cite{Krein:1945}, and Tom{\v c}uk
\cite{Tomcuk:1963}. For a modern treatment of the theory of orthogonal
polynomials on the unit circle and an exhaustive bibliography on
the subject we refer to the forthcoming monumental two-volume treatise  
by Simon \cite{Simon:2004} (see also \cite{Simon:2005}). For fascinating
connections between orthogonal polynomials and random matrix theory we
refer, for instance, to Deift \cite{Deift:1999}.

An extension of \eqref{sb1.1} was developed by Baxter in a
series of papers on Toeplitz forms \cite{Baxter:1960}--\cite{Baxter:1963}
in 1960--63. In these papers the transfer matrix $T$ in \eqref{sb1.3} is
replaced by the more general transfer matrix
\begin{equation}
U(z)=\begin{pmatrix} z & \alpha \\ \beta z & 1 \end{pmatrix} \lb{sb1.9}
\end{equation}
with $\alpha=\alpha(n)$, $\beta=\beta(n)$, subject to the condition
\begin{equation}
\alpha (n)\beta (n)\ne 1\, \text{ for all $n\in\bbN$}. \lb{sb1.10}
\end{equation}
Studying the following extension of \eqref{sb1.1}, 
\begin{equation}
\Psi(z,n) = U(z,n)\Psi(z,n-1), \quad
(z,n)\in\bbT\times\bbN,  \lb{sb1.10a}
\end{equation}
Baxter was led to biorthogonal polynomials on the unit circle with respect
to a complex-valued measure. In this paper we will primarily be concerned
with Baxter's extension \eqref{sb1.10a} of \eqref{sb1.1}.

To simplify our notation in the following, shifts on the lattice are denoted
using superscripts, that is, we write for complex-valued sequences $f$,
\begin{equation}
(S^\pm f)(n)=f^\pm (n)=f(n\pm1), \quad n\in\bbZ, \,\,
\{f(n)\}_{n\in\bbZ}\subset\bbC \lb{sb1.11}
\end{equation}
and apply the analogous convention to $2\times 2$ matrices and their
entries.

In the mid seventies, Ablowitz and Ladik, in a series of papers
\cite{AblowitzLadik:1975}, \cite{AblowitzLadik:1976},
\cite{AblowitzLadik1:1976}, \cite{AblowitzLadik2:1976} (see also
\cite{Ablowitz:1977}, \cite[Sect.\ 3.2.2]{AblowitzClarkson:1991}, 
\cite[Ch.\ 3]{AblowitzPrinariTrubatch:2004}), used inverse scattering
methods to analyze certain integrable differential-difference systems. 
One of their integrable variants of such systems, a discretization of the
AKNS-ZS system, is of the type
\begin{align}
-i\alpha_t-(\alpha^+-2\alpha+\alpha^-)
+\alpha \beta (\alpha^+ +\alpha^-)&=0, \lb{sb1.12} \\
-i\beta_t +(\beta^+-2\beta +\beta^-) -\alpha\beta
(\beta^+ +\beta^-)&=0 \lb{sb1.13}
\end{align}
with $\alpha=\alpha(n,t)$, $\beta=\beta(n,t)$. 
In particular, Ablowitz and Ladik \cite{AblowitzLadik:1976} (see also
\cite[Ch.\ 3]{AblowitzPrinariTrubatch:2004}) showed that in the
focusing case, where $\beta = -\ol\alpha$, and in the defocusing case,
where $\beta = \ol\alpha$ (cf.\ \eqref{sb1.3}),
\eqref{sb1.12} and \eqref{sb1.13} yield the discrete analog of the
nonlinear Schr\"odinger equation
\begin{equation}
-i\alpha_t +2\alpha-(1 \pm |\alpha|^2)(\alpha^+ + \alpha^-)=0. \lb{sb1.17a}
\end{equation}

Algebro-geometric solutions of the AL system \eqref{sb1.12}, \eqref{sb1.13}
have been studied by Ahmad and Chowdury \cite{AhmadChowdhury:1987},
Bogolyubov, Prikarpatskii, and Samoilenko
\cite{BogolyubovPrikarpatskiiSamoilenko:1981},  Bogolyubov and
Prikarpatskii \cite{BogolyubovPrikarpatskii:1982}, Geng, Dai, and Cewen
\cite{GengDaiCao:2003}, Vekslerchik \cite{Vekslerchik:1999}, and especially,
by Miller, Ercolani, Krichever, and Levermore
\cite{MillerErcolaniKricheverLevermore:1995} in an effort to analyze models
describing oscillations in non-linear dispersive wave systems. In 
\cite{MillerErcolaniKricheverLevermore:1995} the authors use the fact that
the AL system \eqref{sb1.12}, \eqref{sb1.13} arises as the compatibility
requirements of the equations
\begin{equation}
\Phi= U \Phi^-,  \quad \Phi_t^-= W \Phi^-. \lb{sb1.20}
\end{equation}
Here $U$ is precisely Baxter's matrix in \eqref{sb1.9} and $W$ is defined
as follows,
\begin{equation}
U(z)=\begin{pmatrix} z & \alpha \\ \beta z & 1 \end{pmatrix},
\quad
W(z)=i\begin{pmatrix}
z-1-\alpha\beta^- & \alpha - \alpha^- z^{-1}\\ \beta^- z -\beta &
1+ \alpha^- \beta -z^{-1} \end{pmatrix}. \lb{sb1.22}
\end{equation}
Thus, the AL system \eqref{sb1.12}, \eqref{sb1.13} is equivalent to the
zero-curvature equations
\begin{equation}
U_t+UW -W^+ U=0.  \lb{sb1.21}
\end{equation}
Miller, Ercolani, Krichever, and Levermore
\cite{MillerErcolaniKricheverLevermore:1995} then performed a thorough
analysis of the solutions $\Phi=\Phi(z,n,t)$ associated with the pair
$(U,W)$ and derived the theta function representations of $\alpha,
\beta$ satisfying the AL system \eqref{sb1.12}, \eqref{sb1.13}. In the
particular focusing and defocusing cases they also discussed periodic
and quasi-periodic solutions $\alpha$ with respect to $n$ and $t$. 

Unaware of the paper \cite{MillerErcolaniKricheverLevermore:1995},
Geronimo and Johnson \cite{GeronimoJohnson:1998} studied the
defocusing case \eqref{sb1.1} in the case where the coefficients $\alpha$
are random variables. They provide a detailed study of the corresponding
Weyl--Titchmarsh functions,
$m_{\pm}$, which  satisfy the Riccati-type equation (for
$z\in\bbC\backslash\bbT$,
$n\in\bbZ$),
\begin{equation}
\alpha(n)m_{\pm}(z,n)m_{\pm}(z,n-1)-m_{\pm}(z,n-1)+zm_{\pm}(z,n)=z\ol
\alpha(n) \lb{sb1.26} 
\end{equation}
(which should be compared to the identical equation \eqref{sb3.13} for the
fundamental function $\phi$ in the defocusing case $\beta=\ol\alpha$). 
These functions take on the values
$|m_+(z)|<1$ and $|m_-(z)|>1$ for $|z|<1$ (cf.\ \cite{GeronimoJohnson:1996},
\cite{GeronimoJohnson:1998}). Utilizing the fact that
$m_{+}$ is a Schur function (i.e., analytic in the open unit disc with
modulus less than one) and the close relation between such functions and the
orthogonality measure $d\sigma_+$, they perform the transformation
\begin{equation}
\hatt\Phi = A\Phi, \quad
A=\f{1}{\sqrt{2}}\begin{pmatrix}
1 & -1\\ i & i\end{pmatrix}. \lb{sb1.27}
\end{equation}
With this change of variables $m_{\pm}$ transform into
\begin{equation}
\hatt m_{\pm}(z,n)=i\frac{1+m_{\pm}(z,n)}{1-m_{\pm}(z,n)}, \quad
z\in\bbC\backslash\bbT, \; n\in\bbZ. \lb{sb1.29}
\end{equation}
The Schur property of $m_{+}$ (equivalently, the relation between Schur
functions,  Caratheodory functions, and positive measures on the
unit circle \cite{Akhiezer:1965}, \cite{Simon:2004a}--\cite{Simon:2005})
implies the standard representation,
\begin{equation}
\hatt m_+(z,n) = i\int^{2\pi}_{0} d\sigma_+(e^{i\theta},n) \, 
\frac{e^{i\theta}+z}{e^{i\theta}-z}, \quad
z\in\bbC\backslash\bbT, \; n\in\bbZ. \lb{sb1.30}
\end{equation}
Under appropriate ergodicity assumptions on $\alpha$ and the
hypothesis of a vanishing Lyapunov exponent on the prescribed spectral
arcs on the unit circle $\bbT$, Geronimo and Johnson
\cite{GeronimoJohnson:1998} showed that the
$m$-functions associated with \eqref{sb1.1} are reflectionless, that is,
$\hatt m_+$ is the analytic continuation of $\hatt m_-$ through the
spectral arcs and vice versa, or equivalently, $\hatt m_\pm$ are the two
branches of an analytic function $\hatt m$ on the hyperelliptic Riemann
surface with branch points given by the end points of the spectral arcs
on $\bbT$. They developed the corresponding spectral theory associated
with \eqref{sb1.1} and the unitary operator it generates in $\ell^2(\bbZ)$
(cf.\ \cite{GeronimoTeplyaev:1994}). This can be viewed as analogous to the
case of real-valued finite-gap potentials for Schr\"odinger operators on
$\bbR$ (cf., e.g., \cite{BelokolosBobenkoEnolskiiItsMatveev:1994},
\cite{GesztesyHolden:2003}) and self-adjoint Jacobi operators on $\bbZ$
(cf., e.g., \cite{BullaGesztesyHoldenTeschl:1997}). In particular, Geronimo
and Johnson \cite{GeronimoJohnson:1998} prove the quasi-periodicity of the
coefficients $\alpha$ in the defocusing case $\beta=\ol\alpha$. Connections
with aspects of integrability, a zero-curvature or Lax formalism, and the
theta function representation of $\alpha$, are not discussed in
\cite{GeronimoJohnson:1998}. The whole topic has been reconsidered in great
detail and partially simplified in the upcoming two-volume monograph by
Simon \cite[Ch.\ 11]{Simon:2004} and aspects
of integrability (Lax pairs, etc.) in the periodic defocusing case will
further be explored by Nenciu and Simon \cite{NenciuSimon:2005}.

The principal contribution of this paper to this circle of ideas is a short 
derivation of theta function formulas for algebro-geometric coefficients
$\alpha, \beta$ associated with Baxter's finite difference
system \eqref{sb1.10a}. Rather than considering solutions of a particular AL
flow such as\eqref{sb1.12}, \eqref{sb1.13}, we will focus on a derivation of
the coupled system of nonlinear difference equations satisfied by
algebro-geometric solutions $\alpha, \beta$ of \eqref{sb1.10a} (a new
result) and its algebro-geometric solutions. In this sense our contribution
represents the analog of determining algebro-geometric coefficients
(generally, complex-valued) in one-dimensional Schr\"odinger and Jacobi
operators and deriving the corresponding Its--Matveev-type theta function
formulas. As a by-product in the special defocusing case $\beta=\ol\alpha$
with $|\alpha(n)|<1$, $n\in\bbZ$, we recover the original result of
Geronimo and Johnson \cite{GeronimoJohnson:1998} that $\alpha$ is
quasi-periodic without the use of Fay's generalized Jacobi variety, double
covers, etc.

In Section \ref{sbs2} we describe our zero-curvature formalism and the
ensuing hierarchy of nonlinear difference equations for $\alpha, \beta$.
Our principal Section \ref{sbs3} then is devoted to a detailed derivation
of the theta function formulas of all algebro-geometric quantities
involved. Appendix \ref{A} collects relevant material on hyperelliptic
curves and their theta functions and introduces the terminology freely used
in Section \ref{sbs3}. 

\section{Zero-Curvature Equations and Hyperelliptic Curves}
\label{sbs2}

In this section we introduce the basic zero-curvature setup for
algebro-geometric solutions of \eqref{sb1.10a}. We follow the approach
employed in \cite{BullaGesztesyHoldenTeschl:1997}, 
\cite{GesztesyHolden:2003}--\cite{GesztesyRatnaseelanTeschl:1996} in the
analogous cases of stationary KdV, AKNS, and Toda solutions.

We start by introducing the complex-valued sequences
\begin{equation}
\{\alpha(n)\}_{n\in\bbZ}, \{\beta(n)\}_{n\in\bbZ}\subset\bbC, \lb{sb2.1}
\end{equation}
and define the recursion relations
\begin{align}
f_{0}&=-2\alpha^+, \quad g_{0}=1, \quad
h_{0}=2\beta, \lb{sb2.22} \\
g_{\ell+1}-g_{\ell+1}^-&=\alpha h_{\ell}^-
         +\beta f_{\ell}, \quad \ell \in \bbN_0,  \lb{sb2.17} \\
f_{\ell+1}^-&=f_\ell-\alpha(g_{\ell+1}+g^-_{\ell+1}),
\quad \ell \in \bbN_0, \lb{sb2.18} \\
h_{\ell+1}&=h_\ell^- +\beta(g_{\ell+1}+g^-_{\ell+1}),
\quad \ell \in \bbN_0,  \lb{sb2.19} 
\end{align}
Here shifts on the lattice are denoted using superscripts as introduced in
\eqref{sb1.11}.

In addition we get the relations
\begin{equation}
g_{\ell+1} - g_{\ell+1}^-=\alpha h_{\ell+1}+\beta f^-_{\ell+1}, \quad 
\ell\in\bbN_0, \lb{sb2.20}
\end{equation} 
which are derived as follows,
\begin{align}
\alpha h_{\ell+1}+\beta f_{\ell+1}^-
&=\alpha h_{\ell}^- +\alpha\beta (g_{\ell+1}+g_{\ell+1}^-)
+\beta f_{\ell} - \alpha\beta (g_{\ell+1}+g_{\ell+1}^-)\no \\
&=\alpha h_{\ell}^- +\beta f_{\ell} 
=g_{\ell+1}-g_{\ell+1}^-, \quad \ell\in\bbN_0, \lb{sb2.16cc}
\end{align}
using relations \eqref{sb2.18}, \eqref{sb2.19}, and \eqref{sb2.17}.

\begin{remark}\lb{sbrem2.1}
One can compute the sequences $\{f_\ell\}$, $\{g_\ell\}$, and
$\{h_\ell\}$ recursively  as follows. Assume that $f_\ell$,
$g_\ell$, and $h_\ell$ are known.  Equation \eqref{sb2.17} is a first
order difference equation in $g_{\ell+1}$ that can be solved directly
and yields a local lattice function. The coefficient  $g_{\ell+1}$ is
determined up to a new constant denoted by $c_{\ell+1}\in\bbC$. Relations
\eqref{sb2.18} and \eqref{sb2.19} then determine $f_{\ell+1}$ and
$h_{\ell+1}$, etc.
       \end{remark}
Explicitly, one obtains
\begin{align}
f_{0}&=-2\alpha^+, \quad 
f_{1}=2\big((\alpha^+)^2\beta+\alpha^+\alpha^{++}\beta^+ -\alpha^{++}\big)
+c_1(-2\alpha^+),\no \\
g_{0}&=1, \quad 
g_{1}=-2\alpha^+\beta+c_1,\lb{sb2.26} \\
h_{0}&=2\beta, \quad 
h_{1}=2\big(-\alpha^+\beta^2-\alpha\beta^-\beta
+\beta^-\big)+ c_1 2\beta, \text{   etc.,} \no
\end{align}
where $\{c_\ell\}_{\ell\in\bbN}\subset\bbC$ denote certain summation
constants. 

Next, assuming $z\in\bbC$, we introduce  the $2\times2$ matrix $U(z)$ by
\begin{equation}
          U(z,n)=\begin{pmatrix}
          z & \alpha(n) \\
           z\beta(n) &1\end{pmatrix},  \quad n\in\bbZ.\lb{sb2.2}
\end{equation}
In addition, we introduce for each fixed $p\in\bbN$ the following $2\times 2$
matrix $V_{p+1}(z)$,
\begin{equation}
V_{p+1}(z,n)=\begin{pmatrix}G^-_{p+1}(z,n) &-F^-_{p}(z,n) \\
H^-_{p+1}(z,n)& -G^-_{p+1}(z,n)\end{pmatrix}, \quad n\in\bbZ, \lb{sb2.3}
\end{equation}
supposing $F_{p}(\dott,n)$ and $G_{p+1}(\dott,n)$, $H_{p+1}(\dott,n)$
to be polynomials of degree $p$ and $p+1$, respectively (cf., however,
Remark \ref{sbr3.0}), with respect to the spectral parameter $z\in\bbC$.

Postulating the stationary zero-curvature condition 
\begin{equation}
U(z,n) V_{p+1}(z,n)-V^+_{p+1}(z,n)U(z,n)=0, \quad p\in\bbN_0, \lb{sb2.6}
\end{equation}
then yields the following fundamental relationships between the
polynomials $F_{p}$, $G_{p+1}$, and $H_{p+1}$,
\begin{align}
F_{p}-zF_{p}^- -\alpha\big(G_{p+1}+G_{p+1}^-\big)&=0,
\lb{sb2.8} \\
z\beta\big(G_{p+1} +G_{p+1}^-\big)+H_{p+1}^- -zH_{p+1}&=0,
\lb{sb2.9} \\
z\big(G_{p+1}^- -G_{p+1}\big)
+\alpha H_{p+1}^- +z\beta F_{p}&=0, \lb{sb2.7} \\
G_{p+1}-G_{p+1}^- -\alpha H_{p+1}-z\beta F_{p}^-&=0.\lb{sb2.10}
\end{align}
Moreover, using relations \eqref{sb2.8}--\eqref{sb2.10} one shows that the
quantity $G_{p+1}^2-F_{p}H_{p+1}$ is a lattice constant and hence the
expression
\begin{equation}
G_{p+1}(z,n)^2-F_{p}(z,n)H_{p+1}(z,n)=R_{2p+2}(z) \lb{sb2.13}
\end{equation}
is an $n$-independent polynomial  of degree $2p+2$ with respect to $z$.
(That $G_{p+1}^2-F_{p}H_{p+1}$, $z\neq 0$, is a lattice constant also
immediately follows from \eqref{sb2.6} taking determinants.)

In order to make the connection between the zero-curvature
formalism and the recursion relation
\eqref{sb2.22}--\eqref{sb2.19}, we now introduce the polynomial ansatz
with respect to the spectral parameter $z$,
\begin{align}
          F_{p}(z)=\sum_{\ell=0}^{p} f_{p-\ell} z^\ell, \quad
          G_{p+1}(z)=\sum_{\ell=0}^{p+1} g_{p+1-\ell} z^\ell, \quad
          H_{p+1}(z)=\sum_{\ell=0}^{p+1} h_{p+1-\ell} z^\ell. \lb{sb2.14}
\end{align}
The stationary zero-curvature condition \eqref{sb2.6} imposes further
restrictions on the coefficients of $V_{p+1}$ that we will now explore.
Since $g_{0}=1$, the quantity
$R_{2p+2}$ in \eqref{sb2.13} is a monic polynomial of degree $2p+2$, that
is,
\begin{equation}
R_{2p+2}(z)=\prod_{m=0}^{2p+1}(z-E_{m}), \quad
\{E_{m}\}_{m=0,\dots,2p+1}\subset\bbC. \lb{sb2.16}
\end{equation}
Next we assume $p\in\bbN$ to avoid cumbersome case distinctions
concerning the trivial case $p=0$. Insertion of \eqref{sb2.14} into
\eqref{sb2.8}--\eqref{sb2.10} then  yields the  relations
\eqref{sb2.22} (normalizing $g_0=1$) and the recursion relations
\eqref{sb2.17}, \eqref{sb2.18}, and \eqref{sb2.19} for $\ell=0,\dots,p-1$.
In addition, one obtains the equations 
\begin{align}
f_{p}-\alpha(g_{p+1}+g_{p+1}^-)&=0, \lb{sb2.16a} \\
\beta(g_{p+1}+g^-_{p+1})+h^-_p -h_{p+1}&=0, \lb{sb2.16A} \\
h_{p+1}^-&=0, \lb{sb2.16b} \\
g_{p+1}^- -g_{p+1}+\alpha h_p^- +\beta f_p&=0, \lb{sb2.16da} \\
\alpha h_{p+1}+g_{p+1}^- -g_{p+1} &=0. \lb{sb2.16c}
\end{align}
Moreover, one infers the relations (cf.\ \eqref{sb2.20})
\begin{align}
g_{\ell}-g_{\ell}^-=\alpha h_\ell+\beta f_\ell^-, \quad 
\ell=0,\dots,p. \lb{sb2.16dd}
\end{align}

Combining \eqref{sb2.16b} and \eqref{sb2.16c}, we first conclude that
$g_{p+1}$ is a lattice constant, that is,
\begin{equation}
g_{p+1}=g_{p+1}^+\in\bbC.  \lb{sb2.16d}
\end{equation}
In addition, using \eqref{sb2.16A}, \eqref{sb2.16b}, and \eqref{sb2.16d} 
one obtains 
\begin{equation}
0=h_{p+1}=h_{p}^-+\beta(g_{p+1}+g_{p+1}^-)=h_{p}^-+2 g_{p+1}\beta,
\lb{sb2.23}
\end{equation}
and hence, $h_p=-2g_{p+1}\beta^+$. Equations \eqref{sb2.16a} and
\eqref{sb2.16d} also yield $f_{p} =2 g_{p+1}\alpha$ in agreement with
\eqref{sb2.16da}. Moreover, \eqref{sb2.16d} is 
consistent with taking $z=0$ in \eqref{sb2.13} which yields
\begin{equation}
g_{p+1}^2=\prod_{m=0}^{2p+1} E_{m}. \lb{sb2.24}
\end{equation}

Thus, the stationary zero-curvature condition
\eqref{sb2.6} is equivalent to a coupled system of nonlinear
difference equations for $\alpha$ and $\beta$ which we write as
\begin{equation}
\sSB_{p}(\alpha,\beta)
=\begin{pmatrix}f_{p}(\alpha,\beta) -2 g_{p+1}\alpha \\
                 h_{p}^-(\alpha,\beta)+2 g_{p+1}\beta \end{pmatrix}=0, 
\quad g_{p+1}=g_{p+1}^+, \lb{sb2.29}
\end{equation}
in honor of the pioneering work by Szeg{\H o} and Baxter in connection
with the transfer matrices \eqref{sb1.3} and \eqref{sb1.9}.  Varying
$p\in\bbN_0$ in
\eqref{sb2.29} then defines the corresponding stationary SB hierarchy of
nonlinear difference equations.  The first few  equations explicitly read
\begin{align}
\sSB_{0}(\alpha,\beta)&=2\begin{pmatrix}-\alpha^+ -g_1\alpha \\
\beta^- +g_1\beta\end{pmatrix}=0, \quad g_1=g_1^+, \no \\
\sSB_{1}(\alpha,\beta)&=2\begin{pmatrix}\alpha^+\alpha^{++}\beta^+
+(\alpha^+)^2\beta
-\alpha^{++}   -c_1\alpha^+ -g_2\alpha \\
-\alpha^-\beta^{--}\beta^- -\alpha(\beta^-)^2+\beta^{--}
       +c_1\beta^- +g_2\beta\end{pmatrix}=0, \lb{sb2.30} \\
&\hspace*{.5cm} g_2=g_2^+, \text{ etc.} \no
\end{align}
By definition, the set of solutions of \eqref{sb2.29}, with $p$ ranging
in $\bbN_0$, represents the class of algebro-geometric solutions associated
with Baxter's finite difference system \eqref{sb1.10a}. The hierarchy of
coupled nonlinear difference equations \eqref{sb2.29} is new.

\begin{remark} \lb{r2.2}
$(i)$ The scaling behavior
$f_\ell(A\alpha,A^{-1}\beta)=Af_\ell(\alpha,\beta)$, 
$g_\ell(A\alpha,A^{-1}\beta)=g_\ell(\alpha,\beta)$,
$h_\ell(A\alpha,A^{-1}\beta)=A^{-1}h_\ell(\alpha,\beta)$, $\ell\in\bbN_0$,
$A\in\bbC\backslash\{0\}$, shows that the stationary SB hierarchy
\eqref{sb2.29} has the scaling invariance,
\begin{equation}
(\alpha,\beta) \mapsto (A\alpha,A^{-1}\beta), \quad
A\in\bbC\backslash\{0\}. \lb{2.31}
\end{equation}
In the special focusing and defocusing cases, where $\beta=-\ol \alpha$
and $\beta=\ol\alpha$, respectively, the scaling constant $A$ in
\eqref{2.31} is further restricted to 
\begin{equation}
|A|=1. \lb{2.32}
\end{equation}
$(ii)$ In the defocusing case $\beta=\ol\alpha$, the compatibility
requirement of the two equations in \eqref{sb2.29} requires the constraint
$|g_{p+1}|^2=1$ and additional spectral theoretic considerations in
connection with the trigonometric moment problem, assuming $|\alpha(n)|<1$,
$n\in\bbZ$, enforce $\{E_m\}_{m=0,\dots,2p+1}\subset\bbT$. The additional
condition of periodicity of $\alpha$ then implies further constraints on 
$\{E_m\}_{m=0,\dots,2p+1}$ (cf.\ \cite[Ch.\ 11]{Simon:2004}). The special
case of  real-valuedness of $\alpha$ also  enforces additional constraints
on $\{E_m\}_{m=0,\dots,2p+1}$. 
\end{remark}

\section{Theta Function Representations} \lb{sbs3}

In this our principal section, we present a detailed study of  
algebro-geometric solutions associated with \eqref{sb1.10a} with special
emphasis on theta function representations of $\alpha, \beta$ and related
quantities. We employ the techniques discussed in
\cite{BullaGesztesyHoldenTeschl:1997} and
\cite{GesztesyHolden:2003} in connection with other integrable systems such
as the KdV, AKNS, and Toda hierarchies.

Throughout this section we suppose
\begin{equation}
\{\alpha(n)\}_{n\in\bbZ}, \{\beta(n)\}_{n\in\bbZ}\subset\bbC,
\quad \alpha(n)\beta(n)\neq 0,1, \; n\in\bbZ,  \lb{sb3.0}
\end{equation}
and assume \eqref{sb2.22}--\eqref{sb2.19}, \eqref{sb2.6},
\eqref{sb2.14}.  Moreover, we freely employ the formalism developed in
Section \ref{sbs2}, keeping $p\in\bbN_0$ fixed.

Returning to \eqref{sb2.16} we now introduce the hyperelliptic curve
$\calK_{p}$ with nonsingular affine part defined by
\begin{align}
          &\calK_{p}\colon \calF_{p}(z,y)=y^2-R_{2p+2}(z)=0,  \lb{sb3.1} \\
          &R_{2p+2}(z)=\prod_{m=0}^{2p+1} (z-E_{m}), \quad
          \{E_{m}\}_{m=0.\dots,2p+1}\subset\bbC\backslash\{0\},
\lb{sb3.2}
\\ &E_{m}\neq E_{m'} \text{ for } m\neq m', \,\, m,m'=0,\dots,2p+1.
\lb{sb3.2a}
\end{align}
Equations \eqref{sb3.0}--\eqref{sb3.2a} are assumed for the remainder of
this section. We compactify $\calK_{p}$ by adding two points
$\Pinfp$ and $\Pinfm$, $\Pinfp\neq\Pinfm$, at infinity, still denoting its
projective closure by $\calK_{p}$.  Finite points $P$ on
$\calK_{p}$ are denoted by $P=(z,y)$ where $y(P)$ denotes the
meromorphic function on $\calK_{p}$ satisfying
$\calF_{p}(z,y)=0$.  The complex structure on $\calK_{p}$ is then
defined in a standard manner and $\calK_{p}$ has topological genus $p$.
Moreover, we use the involution
\begin{equation}
          *\colon \calK_{p}\to \calK_{p}, \quad P=(z,y)\mapsto
          P^{*}=(z,-y),\; \Pinfpm\mapsto\Pinfpm^*=\Pinfmp. \lb{sb3.3}
\end{equation}
For further properties and notation concerning hyperelliptic curves
we refer to Appendix \ref{A}.

\begin{remark} \lb{sbr3.0}
The assumption $\alpha(n)\neq 0$, $\beta(n)\neq 0$, $n\in\bbZ$, in
\eqref{sb3.0} is not an essential one. It has the advantage of
guaranteeing that for all $n\in\bbZ$, $F_p(\dott,n)$ and
$H_{p+1}(\dott,n)$ are polynomials of degree $p$ and $p+1$,
respectively. If
$\alpha^+(n_0)=0$ (resp., $\beta(n_0)=0$) for some $n_0\in\bbZ$, then
$F_p(\dott,n_0)$ has at most degree $p-1$ (resp.,
$H_{p+1}(\dott,n_0)$ has at most degree $p$). The latter
$n$-dependence of the degree of the polynomials $F_p$ and
$H_{p+1}$ enforces numerous case distinctions in connection with
our fundamental function $\phi$ in \eqref{sb3.7} below. For simplicity we
will in almost all situations avoid these cumbersome case distinctions and
hence assume $\alpha(n)\neq 0$,
$\beta(n)\neq 0$, $n\in\bbZ$ throughout this section. (For an exception
see Remark \ref{sbr3.5A}.) In the extreme case that $\alpha\equiv 0$ (i.e.,
$\alpha(n)=0$ for all $n\in\bbZ$), then $F_p\equiv 0$ and hence the curve
$\calK_p$ becomes singular as $R_{2p+2}(z)=G_{p+1}(z,n)^2$, $z\in\bbC$,  
by \eqref{sb2.13}, and thus the branch points of $\calK_{p}$ necessarily
occur in pairs. (In addition, $G_{p+1}(z,n)$ then becomes independent of
$n\in\bbZ$.) The same argument applies to $\beta\equiv 0$ since then
$H_{p+1}\equiv 0$. For this reason the trivial cases $\alpha\equiv 0$ and
$\beta\equiv 0$ in \eqref{sb3.0} are excluded in the remainder of this
paper. Finally, in order to avoid numerous case distinctions in connection
with the trivial case $p=0$, we shall assume $p\ge 1$ for the remainder of
this section (with the exception of Example \ref{sbexample3.9}).
\end{remark}

In the following, the zeros of the polynomials $F_{p}(\dott,n)$
and $H_{p+1}(\dott,n)$ (cf.\ \eqref{sb2.14}) will play a special
role. We denote them by $\{\mu_j(n)\}_{j=1,\dots,p}$ and
$\{\nu_\ell(n)\}_{\ell=0,\dots,p}$ and hence write
\begin{equation}
          F_{p}(z)=-2\alpha^+\prod_{j=1}^p(z-\mu_{j}), \quad
          H_{p+1}(z)=2\beta\prod_{\ell=0}^{p}(z-\nu_{\ell}). \lb{sb3.4}
\end{equation}
In addition, we lift these zeros to $\calK_p$ by introducing 
\begin{align}
          \hat\mu_{j}(n)=(\mu_{j}(n),
          G_{p+1}(\mu_{j}(n),n))\in\calK_{p},\quad
          j=1,\dots,p, \; n\in\bbZ,  \lb{sb3.5} \\
           \hat\nu_{\ell}(n)=(\nu_{\ell}(n),-
G_{p+1}(\nu_{\ell}(n),n))\in\calK_{p}, \quad
          \ell=0,\dots,p, \; n\in\bbZ. \lb{sb3.6}
\end{align}
We recall that  $h_{p+1}=0$ (cf.\  \eqref{sb2.16b}). Hence  we may choose
\begin{equation}
\nu_0(n)=0, \quad n\in\bbZ. \lb{sb3.6aa}
\end{equation}
Define
\begin{equation}
\Pzpm=(0,\pm G_{p+1}(0))=(0,\pm g_{p+1}), \lb{sb3.6ab}
\end{equation}
where
\begin{equation}
y(\Pzpm)=\pm g_{p+1}, \quad g_{p+1}^2=\prod_{m=0}^{2p+1}E_m. \lb{sb3.6aB}
\end{equation}
We emphasize that $\Pzpm$ and $\Pinfpm$ are not necessarily on the same sheet
of $\calK_p$. The actual sheet on which $\Pzpm$ lie depends on the sign of
$g_{p+1}$. Thus, one obtains
\begin{equation}
\hnu_0=\Pzm. \lb{sb3.6ac}
\end{equation}
The branch of $y(\dott)$ near $\Pinfpm$ is fixed
according to
\begin{equation}
\lim_{\substack{\abs{z(P)}\to\infty\\P\to\Pinfpm}}\f{y(P)}{
G_{p+1}(z(P))}=\lim_{\substack{\abs{z(P)}\to\infty\\P\to\Pinfpm}}
\f{y(P)}{z(P)^{p+1}}=\mp 1. \lb{sb3.45}
\end{equation}

Next, we introduce the fundamental meromorphic function $\phi$ on
$\calK_{p}$ by
\begin{equation}
          \phi(P,n)=\frac{y+G_{p+1}(z,n)}{F_{p}(z,n)} 
          =\frac{-H_{p+1}(z,n)}{y-G_{p+1}(z,n)}, \quad
          P=(z,y)\in\calK_{p}, \, n\in\bbZ \lb{sb3.7}
\end{equation}
with divisor $(\phi(\dott,n))$ (cf.\ the notation for divisors
introduced in \eqref{A.17} and \eqref{A.18}) given by
\begin{equation}
          (\phi(\dott,n))=\calD_{\Pzm\hunu(n)}-\calD_{\Pinfm\humu(n)}.
          \lb{sb3.8}
\end{equation}
Here we abbreviated (cf.\ \eqref{A.17}, \eqref{A.18})
\begin{equation}
\humu=\{\hmu_{1},\dots,\hmu_{p}\}, \,
\hunu=\{\hnu_{1},\dots,\hnu_{p}\}\in\sym^{p}\calK_{p}. 
\lb{sb3.9}
\end{equation}

The stationary Baker--Akhiezer vector $\Psi(P,n,n_{0})$ is defined on
$\calK_{p}$ by
\begin{align}
\Psi(P,n,n_{0})&=\begin{pmatrix}
          \psi_{1}(P,n,n_{0})\\ \psi_{2}(P,n,n_{0})\end{pmatrix}, \lb{3.10} \\
\psi_{1}(P,n,n_{0})&=\begin{cases}
\prod_{m=n_{0}+1}^{n}
(z+\alpha(m)\phi^-(P,m)), & n\ge n_{0}+1, \\
1, & n=n_{0}, \\
       \prod_{m=n+1}^{n_{0}}
(z+\alpha(m)\phi^-(P,m))^{-1}, & n\le n_{0}-1,
\end{cases} \lb{sb3.11} \\
\psi_{2}(P,n,n_{0})
&=\phi(P,n_{0})  \begin{cases}
\prod_{m=n_{0}+1}^{n}
(z\beta(m)\phi^-(P,m)^{-1}+1), & n\ge n_{0}+1, \\
1, & n=n_{0}, \\
       \prod_{m=n+1}^{n_{0}}
(z\beta(m)\phi^-(P,m)^{-1} +1)^{-1}, & n\le n_{0}-1,
\end{cases} \no \\
&\hspace*{5.5cm} P\in\calK_{p}, \,\,(n,n_{0})\in\bbZ^2. \lb{sb3.12}
\end{align}
Clearly $\Psi(\dott,n,n_{0})$ is meromorphic on $\calK_{p}$ since
$\phi(\dott,m)$ is meromorphic on $\calK_{p}$.  Fundamental
properties of $\phi$ and $\Psi$ are summarized next.

\begin{lemma}\lb{sblemma3.1}
Suppose $\alpha, \beta \subset\bbC$ satisfy \eqref{sb3.0} and the $p$th
stationary SB system \eqref{sb2.29}. Moreover, assume
\eqref{sb3.1}--\eqref{sb3.2a} and let
$P=(z,y)\in\calK_{p}\backslash\{\Pinfp,\Pinfm\}$,
$(n,n_{0})\in\bbZ^2$. Then $\phi$ satisfies the Riccati-type
equation
\begin{equation}
     \alpha \phi(P)\phi^-(P)-\phi^-(P)+z\phi(P) =z\beta \lb{sb3.13}
\end{equation}
and $\Psi$ fulfills
\begin{align}
&\psi_{2}(P,n,n_{0})/\psi_{1}(P,n,n_{0})=\phi(P,n), \lb{sb3.17} \\
& \Psi(P,n,n_{0})=U(z,n)\Psi^-(P,n,n_{0}), \lb{sb3.18} \\
& -y \Psi^- (P,n,n_{0})=V_{p+1}(z,n)\Psi^-(P,n,n_{0}). \lb{sb3.19}
\end{align}
\end{lemma}
\begin{proof} Using $y^2=G_{p+1}^2-F_{p}H_{p+1}$ 
(cf.\ \eqref{sb2.13}, \eqref{sb3.1}) and \eqref{sb3.7}, the left-hand side
of \eqref{sb3.13} can be rewritten as follows
\begin{align}
&\alpha\phi\phi^- - \phi^- +z\phi-z\beta 
=(F_{p}F_{p}^-)^{-1}\big[\alpha\big(G_{p+1}^2-F_{p}H_{p+1}+y(G_{p+1}
+G_{p+1}^-) \no \\
& \quad + G_{p+1}G_{p+1}^- \big) -(y+G_{p+1}^-)F_{p}
+z(y+G_{p+1})F_{p}^- -z\beta F_{p}F_{p}^- \big]. \lb{sb3.25}
\end{align}
Insertion of \eqref{sb2.8} and \eqref{sb2.10} into \eqref{sb3.25} then
shows that the right-hand side of \eqref{sb3.25} vanishes. This
proves \eqref{sb3.13}. Relation
\eqref{sb3.17} is proven inductively as follows.  Since it holds for
$n=n_{0}$ by
\eqref{sb3.11}, \eqref{sb3.12} we assume that
\begin{equation}
          \psi_{2}(P,m,n_{0})/\psi_{1}(P,m,n_{0})=\phi(P,m), \quad
          m=n_{0},\dots,n-1. \lb{sb3.30}
\end{equation}
Then combining \eqref{sb3.11}, \eqref{sb3.12}, and \eqref{sb3.30},
one obtains
\begin{equation}
          \frac{\psi_{2}(P,n,n_{0})}{\psi_{1}(P,n,n_{0})}=\phi^-(P,n)
          \frac{z\beta(n)\phi^-(P,n)^{-1}+1}
          {z+\alpha(n)\phi^-(P,n)}, \lb{sb3.31}
\end{equation}
and hence
\begin{equation}
\alpha(n)\frac{\psi_{2}(P,n,n_{0})}{\psi_{1}(P,n,n_{0})}\phi^-(P,n)-
\phi^-(P,n) +z\frac{\psi_{2}(P,n,n_{0})}{\psi_{1}(P,n,n_{0})}
-z\beta(n)=0. \lb{sb3.32}
\end{equation}
Comparison with \eqref{sb3.13} (cf.\ also \eqref{sb3.25}) then proves
\eqref{sb3.17} for all $n\ge n_{0}$. The case $n\le n_{0}-1$ is
proven analogously.  By \eqref{sb3.11} and \eqref{sb3.12} one infers
\begin{align}
\psi_{1}(P,n,n_{0})&=
[z+\alpha(n)\phi^-(P,n)]\psi^-_{1}(P,n,n_{0}) \no \\
&=z\psi^-_{1}(P,n,n_{0})+\alpha(n)\psi^-_{2}(P,n,n_{0}),
       \lb{sb3.33} \\
\psi_{2}(P,n,n_{0})
&= [z\beta(n)\phi^-(P,n)^{-1}+1]\psi^-_{2}(P,n,n_{0})
\no \\
&=z\beta(n)\psi^-_{1}(P,n,n_{0})+\psi^-_{2}(P,n,n_{0}),
          \lb{sb3.34}
\end{align}
by \eqref{sb3.17}.  This proves \eqref{sb3.18}.  An application of
\eqref{sb3.7} implies
\begin{align}
          V_{p+1}\Psi^- = \begin{pmatrix}
          G_{p+1}^- \psi_{1}^- -  F_{p}^-\psi_{2}^- \\
          H_{p+1}^- \psi_{1}^- -  G_{p+1}^-\psi_{2}^-
\end{pmatrix} 
=\begin{pmatrix}(G_{p+1}^- -F_{p}^- \phi^-)\psi_{1}^- \\
(H_{p+1}^{-}(\phi^-)^{-1}-G_{p+1}^-) \psi_{2}^-
       \end{pmatrix}=-y \Psi^- \lb{sb3.35}
\end{align}
and hence \eqref{sb3.19}.  
\end{proof}

We note that the Riccati-type equation \eqref{sb3.13} for $\phi$ coincides
with that of $m_\pm$ in \eqref{sb1.26} in the defocusing case
$\beta=\ol\alpha$.

Next, we derive trace formulas for $\alpha$ and $\beta$ in terms of the
zeros $\mu_{j}$ and $\nu_{j}$ of $F_p$ and $H_{p+1}$, respectively. For
simplicity we just record the simplest case below.

\begin{lemma}\lb{sblemma3.3}
Suppose that $\alpha, \beta \subset\bbC$ satisfy \eqref{sb3.0} and the
$p$th stationary SB system \eqref{sb2.29}. Then,
\begin{equation}
          \frac{\alpha}{\alpha^+}=\frac{(-1)^{p+1}}{g_{p+1}}
          \prod_{j=1}^p \mu_{j}, \quad 
       \frac{\beta^+}{\beta}=\frac{(-1)^{p+1}}{g_{p+1}}
          \prod_{\ell=1}^p \nu_{\ell}.\lb{sb3.40}
\end{equation}
\end{lemma}
\begin{proof}
Combining \eqref{sb2.14}, $f_p=2g_{p+1}\alpha$, and \eqref{sb3.4} yields
\begin{equation}
2g_{p+1}\alpha=f_{p}=f_{0}(-1)^p\prod_{j=1}^p \mu_{j}=
-2\alpha^+(-1)^p\prod_{j=1}^p\mu_{j}. \lb{sb3.42} 
\end{equation}
Using $h_p=-2g_{p+1}\beta^+$, the case $\beta^+/\beta$ is analogous.
\end{proof}

The following result describes the asymptotic behavior of $\phi$
as $P\to\Pinfpm$ and $P\to \Pzpm$.

\begin{lemma}\lb{sblemma3.4}
Suppose that $\alpha, \beta \subset\bbC$ satisfy \eqref{sb3.0} and the
$p$th SB system \eqref{sb2.29}. In addition, assume
\eqref{sb3.1}--\eqref{sb3.2a} and let
$P=(z,y)\in\calK_{p}\backslash\{\Pinfp,\Pinfm\}$ and
$(n,n_0)\in\bbZ^2$. Then,
\begin{align}
\phi(P)&\underset{\zeta\to 0}{=} \begin{cases}
\beta+(1-\alpha\beta)\beta^- \zeta+\Oh(\zeta^{2}) &
\text{as $P\to\Pinfp$}, \\
-(\alpha^+)^{-1}\zeta^{-1}
+(1-\alpha^+\beta^+)\alpha^{++}(\alpha^+)^{-2} +\Oh(\zeta)
& \text{as $P\to\Pinfm$}, \end{cases} \no \\
& \hspace*{8.6cm} \zeta=1/z, \lb{sb3.46} \\
\phi(P)&\underset{\zeta\to 0}{=} \begin{cases}
(\alpha)^{-1}- (1-\alpha\beta)\alpha^- (\alpha)^{-2} \zeta
+\Oh(\zeta^2) & \text{as $P\to\Pzp$}, \\
-\beta^+ \zeta
-(1-\alpha^+\beta^+)\beta^{++}\zeta^2+\Oh(\zeta^{3}) &
\text{as $P\to\Pzm$}, \end{cases} \quad \zeta=z. \lb{sb3.47}
\end{align}
\end{lemma}
\begin{proof}  Inserting the ansatz
\begin{equation}
\phi(P,n)=\begin{cases} 
\phi_{-1}(n)z+\phi_{0}(n)+\phi_{1}(n)z^{-1}+\Oh(z^{-2}) & 
\text{as $z\to\infty$},  \\
\phi_{0}(n)+\phi_{1}(n)z+\Oh(z^{2}) & \text{as $z\to 0$} 
\end{cases} \lb{sb3.51}
\end{equation}
into the Riccati-type equation \eqref{sb3.13} produces \eqref{sb3.46} and
\eqref{sb3.47}.
\end{proof}

Using \eqref{sb3.7}--\eqref{sb3.12} one obtains for the divisor
$(\psi_{j}(\dott,n,n_{0}))$ of the meromorphic functions
$\psi_{j}(\dott,n,n_{0})$, $j=1,2$, 
\begin{align}
(\psi_{1}(\dott,n,n_{0}))&=\calD_{\humu(n)}
-\calD_{\humu(n_{0})}+(n-n_{0})(\calD_{\Pzm}-\calD_{\Pinfp}),
          \lb{sb3.48}  \\
(\psi_{2}(\dott,n,n_{0}))&=\calD_{\Pzm\hunu(n)}
-\calD_{\Pinfm\humu(n_{0})}+(n-n_{0})(\calD_{\Pzm}-\calD_{\Pinfp}).
          \lb{sb3.49}
\end{align}

Next, we briefly consider the asymptotic behavior of
$\phi$ in the case where the conditions $\alpha(n)\beta(n)\neq
0$ are violated for some $n\in\bbZ$.

\begin{remark} \lb{sbr3.5A} 
First we note that if $\alpha^+\neq 0$, then by \eqref{sb3.46} no pole
$\hmu_j$ of $\phi$ hits the point $\Pinfm$. Similarly, by \eqref{sb3.47},
$\Pzm=\hnu_0$ is a zero of $\phi$. The case $\beta(n)=0$ for
some $n\in\bbZ$ poses no difficulty and \eqref{sb3.46} and \eqref{sb3.47}
extend continuously in this case. The case $\alpha(n)=0$ for some
$n\in\bbZ$ is more involved and causes higher order poles in \eqref{sb3.46}
and \eqref{sb3.47}. An explicit calculation yields ($\zeta=1/z$)
\begin{equation}
\phi(P)\underset{\zeta\to 0}{=} \begin{cases}
\Oh(1) & \text{as $P\to\Pinfp$} \\
-(\alpha^{++})^{-1}\zeta^{-2} +\Oh(\zeta^{-1})
& \text{as $P\to\Pinfm$} \end{cases} \; \text{ if $\alpha^+=0$,
$\alpha^{++}\neq 0$}. \lb{sb3.55B}  
\end{equation}
Thus, if $\alpha^+=0$, $\alpha^{++}\neq 0$, one of the poles $\hmu_j$ of
$\phi$ hits the point $\Pinfm$. However, still no pole of $\phi$ hits
$\Pinfp$. Similarly, using
\begin{align}
f_p&=2g_{p+1}\alpha, \quad 
f_{p-1}=2g_{p+1}(\alpha^- -\alpha^2\beta^+-\alpha^-\alpha\beta)
+2C\alpha, \no \\
g_p&=-2g_{p+1}\alpha\beta^+ + g_{p+1}^{-1}2(2p+1)c_1,  \lb{sb3.55C} \\
h_p&=-2g_{p+1}\beta^+, \quad 
h_{p-1}=2g_{p+1}(-\beta^{++} +\alpha\beta\beta^+
+\alpha^-\beta^2)-2C\beta, \no
\end{align}
one derives ($\zeta=z$)
\begin{equation}
\phi(P)\underset{\zeta\to 0}{=} \begin{cases}
(\alpha^-)^{-1}\zeta^{-1} +\Oh(1) & \text{as $P\to\Pzp$} \\
\Oh(\zeta) & \text{as $P\to\Pzm$} \end{cases}
\; \text{ if $\alpha=0$, $\alpha^{-}\neq 0$}. \lb{sb3.55E} 
\end{equation}
Thus, if $\alpha=0$, $\alpha^{-}\neq 0$, one of the poles $\hmu_j$ of
$\phi$ hits the point $\Pzp$. In addition, $\Pzm$ remains a zero of $\phi$.
\end{remark}

Since nonspecial divisors will play a fundamental role in this section, we
now take a closer look at them.

\begin{lemma} \lb{sbl3.4}
Suppose that $\alpha, \beta \subset \bbC$ satisfy \eqref{sb3.0} and the
$p$th stationary SB system \eqref{sb2.29}. Moreover, assume
\eqref{sb3.1}--\eqref{sb3.2a} and let $n\in\bbZ$. Let $\calD_{\humu}$,
$\humu=\{\hmu_1,\dots,\hmu_p\}$ and $\calD_{\hunu}$,
$\hunu=\{\hunu_1,\dots,\hunu_p\}$, be the pole and zero divisors of degree
$p$, respectively, associated with $\alpha, \beta$ and $\phi$ defined
according to
\eqref{sb3.5}, \eqref{sb3.6}, that is,
\begin{align}
\hat\mu_j (n)
&= (\mu_j (n),G_{p+1}(\mu_j (n), n)), \quad j=1,\dots,p, \; n\in\bbZ,
\lb{3.44} \\
\hat\nu_j (n)
&= (\nu_j (n),-G_{p+1}(\nu_j (n), n)), \quad j=1,\dots,p, \; n\in\bbZ.
\lb{3.45}
\end{align}
Then $\calD_{\humu(n)}$ and $\calD_{\hunu(n)}$ are nonspecial for all
$n\in\bbZ$.
\end{lemma}
\begin{proof}
We provide a detailed proof in the case of $\calD_{\humu(n)}$.
By Theorem \ref{thm3}, $\calD_{\humu(n)}$ is special if and only if
$\{\hmu_1(n),\dots,\hmu_p(n)\}$ contains at least one pair of the type
$\{\hat\mu(n),\hat\mu^*(n)\}$. Hence $\calD_{\humu(n)}$ is certainly
nonspecial as long as the projections $\mu_j(n)$ of $\hmu_j(n)$ are
mutually distinct, $\mu_j(n) \neq \mu_k(n)$ for $j\neq k$. On the other
hand, if two or more projections coincide for some $n_0\in\bbZ$, for
instance,
\begin{equation}
\mu_{j_1}(n_0)=\cdots=\mu_{j_N}(n_0)=\mu_0, \quad  N\in\{2,\dots,p\},
\lb{sb3.45AA}
\end{equation}
then $G_{p+1}(\mu_0,n_0)\neq 0$ as long as $\mu_0\notin
\{E_0,\dots,E_{2p+1}\}$. This fact immediately follows from \eqref{sb2.13}
since $F_p(\mu_0,n_0)=0$ but $R_{2p+2}(\mu_0)\neq 0$ by hypothesis. In
particular, $\hmu_{j_1}(n_0),\dots,\hmu_{j_N}(n_0)$ all meet on the same
sheet since
\begin{equation}
\hmu_{j_r}(n_0)=(\mu_0,G_{p+1}(\mu_0,n_0)), \quad r=1,\dots,N
\lb{sb3.45a}
\end{equation}
and hence no special divisor can arise in this manner. It remains to
study the case where two or more projections collide at a branch point,
say at $(E_{m_0},0)$ for some $n_0\in\bbZ$. In this case one concludes
$F_p(z,n_0)\underset{z\to E_{m_0}}{=}\Oh\big((z-E_{m_0})^2\big)$ and
\begin{equation}
G_{p+1}(E_{m_0},n_0)=0 \lb{sb3.45c}
\end{equation}
using again \eqref{sb2.13} and $F_p(E_{m_0},n_0)=R_{2p+2}(E_{m_0})=0$.
Since $G_{p+1}(\dott,n_0)$ is a polynomial (of degree $p+1$),
\eqref{sb3.45c} implies
$G_{p+1}(z,n_0)\underset{z\to E_{m_0}}{=}\Oh((z-E_{m_0}))$. 
Thus, using \eqref{sb2.13} once more, one obtains the contradiction,
\begin{align}
\Oh\big((z-E_{m_0})^2\big)&\underset{z\to E_{m_0}}{=} R_{2p+2}(z)
\lb{sb3.45e}
\\ &\underset{z\to
E_{m_0}}{=} (z-E_{m_0})\bigg(\prod_{\substack{m=1\\m\neq m_0}}^{2p+1}
\big(E_{m_0}-E_m\big) +\Oh(z-E_{m_0})\bigg). \no
\end{align}
Consequently, at most one $\hmu_j(n)$ can hit a branch point at a time
and again no special divisor arises. Finally, by \eqref{sb3.46},
$\hat\mu_j(n)$ never reaches the points $\Pinfp$. Hence if some
$\hmu_j(n)$ tend to infinity, they all necessarily converge to $\Pinfm$.
Again no special divisor can arise in this manner.

The proof for $\calD_{\hunu(n)}$ is completely analogous
(replacing $F_p$ by $H_{p+1}$ and noticing that by \eqref{sb3.46},
$\phi$ has no zeros near $\Pinfpm$), thereby completing the proof.
\end{proof}

\begin{remark} \lb{sbr3.9}
For simplicity we assumed $\alpha(n)\neq 0$, $\beta(n)\neq 0$, $n\in\bbZ$,
in Lemma \ref{sbl3.4}. However, the asymptotic behavior in \eqref{sb3.55B}
(resp., \eqref{sb3.55E}) shows that no special divisors can be created at
infinity (resp., zero) and hence the results of Lemma \ref{sbl3.4} extend by
continuity to the situation considered in Remark \ref{sbr3.5A}. In
particular, it extends to the case where $\beta(n_0)=0$ for some
$n_0\in\bbZ$. The case $\alpha(n_0)=0$ for some $n_0\in\bbZ$ is more
involved and requires more and more case distinctions as is clear from
Remark \ref{sbr3.5A}, but the pattern persists.
\end{remark}

Next we turn to the representation of $\phi$, $\Psi$, $\alpha$, and
$\beta$ in terms of the Riemann theta function associated with
$\calK_{p}$.  We freely use the notation established in Appendix
\ref{A}, assuming $\calK_{p}$ to be nonsingular as in
\eqref{sb3.1}--\eqref{sb3.2a}. To avoid the trivial case $p=0$
(considered separately in Example \ref{sbexample3.9}), we assume
$p\in\bbN$ for the remainder of this argument.

We choose a fixed base point
$Q_{0}\in\calK_{p}\backslash\{\Pinfp,\Pinfm, \Pzp,\Pzm\}$, in
fact, we will  choose a branch point for convenience,
$Q_{0}\in\calB(\calK_{p})$. Moreover we denote by
$\omega^{(3)}_{P_{1},P_{2}}$ a normal differential of the third  kind
(cf.\ \eqref{a36}, \eqref{a37}) with simple poles at
$P_{1}$ and $P_{2}$ with residues $1$ and $-1$, respectively. Explicitly,
one computes for $\omega_{\Pzm, \Pinfm}^{(3)}$ and 
$\omega_{\Pzm,\Pinfp}^{(3)}$ the following expressions
\begin{equation}
\omega_{\Pzm, \Pinfpm}^{(3)}=\f{y+y_{0,-}}{z}\,\f{dz}{2y}
\mp \f{1}{2y}\prod_{j=1}^p (z-\lambda_{\pm,j})dz, \quad
\Pzm=(0, y_{0,-})=(0,-g_{p+1}), \lb{sb3.57a}
\end{equation}
where $\{\lambda_{\pm,j}\}_{j=1,\dots,p}$ are uniquely determined by the
normalization
\begin{equation}
\int_{a_{j}}\omega_{\Pzm, \Pinfpm}^{(3)}=0, \quad
j=1,\dots,p.  \lb{sb3.57b}
\end{equation}
The explicit formula \eqref{sb3.57a} then implies the following
asymptotic expansions (using the local coordinate $\zeta=z$ near $\Pzpm$
and $\zeta=1/z$ near $\Pinfpm$),
\begin{align}
\int_{Q_0}^P \omega_{\Pzm, \Pinfm}^{(3)}&\underset{\zeta\to 0}{=}
\left\{\begin{matrix} 0 \\\ln(\zeta) \end{matrix} \right\}
+\omega_0^{0,\pm}(\Pzm,\Pinfm) + \Oh(\zeta)
\text{  as $P\to \Pzpm$}, \lb{sb3.57c} \\
\int_{Q_0}^P \omega_{\Pzm, \Pinfm}^{(3)}&\underset{\zeta\to 0}{=}
\left\{\begin{matrix} 0 \\-\ln(\zeta) \end{matrix} \right\}
+\omega^{\infty_\pm}_0(\Pzm,\Pinfm) + \Oh(\zeta)
\text{  as $P\to \Pinfpm$}, \lb{sb3.57d} \\
\int_{Q_0}^P \omega_{\Pzm, \Pinfp}^{(3)}&\underset{\zeta\to 0}{=}
\left\{\begin{matrix} 0 \\\ln(\zeta) \end{matrix} \right\}
+\omega_0^{0,\pm}(\Pzm,\Pinfp) + \Oh(\zeta)
\text{  as $P\to \Pzpm$}, \lb{sb3.57e} \\
\int_{Q_0}^P \omega_{\Pzm, \Pinfp}^{(3)}&\underset{\zeta\to 0}{=}
\left\{\begin{matrix} -\ln(\zeta) \\ 0\end{matrix} \right\}
+\omega^{\infty_\pm}_0(\Pzm,\Pinfp) + \Oh(\zeta)
\text{  as $P\to \Pinfpm$}. \lb{sb3.57f}
\end{align}
Here $Q_0\in\calB(\calK_{n})$ is a fixed base point and we
agree to choose the same path of integration from $Q_0$ to $P$ in all
Abelian integrals in this section.

\begin{lemma} \lb{sbl3.8a}
With $\omega_0^{\infty_{\sigma}}(\Pzm,\Pinfpm)$ and
$\omega_0^{0,\sigma'}(\Pzm,\Pinfpm)$, $\sigma, \sigma' \in\{+,-\}$, defined as
in \eqref{sb3.57c}--\eqref{sb3.57f} one has 
\begin{align}
\exp\big[&\omega_0^{0,-}(\Pzm,\Pinfpm)
-\omega_0^{\infty_+}(\Pzm,\Pinfpm) \no\\
& -\omega_0^{\infty_-}(\Pzm,\Pinfpm)
+\omega_0^{0,+}(\Pzm,\Pinfpm)\big]=1. \lb{sb3.57g}
\end{align}
\end{lemma}
\begin{proof}
Pick $Q_{1,\pm}=(z_1,\pm y_1)\in\calK_n\backslash\{\Pinfpm\}$ in a
neighborhood of $\Pinfpm$ and $Q_{2,\pm}=(z_2,\pm
y_2)\in\calK_n\backslash\{\Pzpm\}$ in a neighborhood of $\Pzpm$.
Without loss of generality we may assume that $\Pinfp$ and $\Pzp$ lie
on the same sheet. Then by \eqref{sb3.57a},
\begin{align}
&\int_{Q_0}^{Q_{2,-}} \omega^{(3)}_{\Pzm,\Pinfm}
-\int_{Q_0}^{Q_{1,+}} \omega^{(3)}_{\Pzm,\Pinfm}
-\int_{Q_0}^{Q_{1,-}} \omega^{(3)}_{\Pzm,\Pinfm}
+\int_{Q_0}^{Q_{2,+}} \omega^{(3)}_{\Pzm,\Pinfm} \no \\
& \quad = \int_{Q_0}^{Q_{2,+}} \f{dz}{z}
  -\int_{Q_0}^{Q_{1,+}} \f{dz}{z} =\ln(z_2)-\ln(z_1)+2\pi ik, \lb{sb3.57h}
\end{align}
for some $k\in\bbZ$. On the other hand, by
\eqref{sb3.57c}--\eqref{sb3.57f} one obtains
\begin{align}
&\int_{Q_0}^{Q_{2,-}} \omega^{(3)}_{\Pzm,\Pinfm}
-\int_{Q_0}^{Q_{1,+}} \omega^{(3)}_{\Pzm,\Pinfm}
-\int_{Q_0}^{Q_{1,-}} \omega^{(3)}_{\Pzm,\Pinfm}
+\int_{Q_0}^{Q_{2,+}} \omega^{(3)}_{\Pzm,\Pinfm} \no \\
& \quad = \ln(z_2)+\ln(1/z_1)+\omega_0^{0,-}(\Pzm,\Pinfm)
-\omega_0^{\infty_+}(\Pzm,\Pinfm)
-\omega_0^{\infty_-}(\Pzm,\Pinfm) \no \\
& \qquad +\omega_0^{0,+}(\Pzm,\Pinfm) + \Oh(z_2)+\Oh(1/z_1), \lb{sb3.57i}
\end{align}
and hence the part of \eqref{sb3.57g} concerning
$\omega^{(3)}_{\Pzm,\Pinfm}$ follows. The corresponding result for
$\omega^{(3)}_{\Pzm,\Pinfp}$ is proved analogously.
\end{proof}


In the following it is convenient to use the abbreviation
\begin{align}
\uz(P,\ul Q)=\uxi_{Q_{0}}-\ua_{Q_{0}}(P)+\ual_{Q_{0}}(\calD_{\ul Q}),
\quad P\in\calK_{p}, \; \ul
Q=\{Q_{1},\dots,Q_{p}\}\in\sym^{p}\calK_{p}. \lb{sb3.58}
\end{align}

\begin{theorem}\lb{sbt3.6}
Suppose that $\alpha, \beta \subset\bbC$ satisfy \eqref{sb3.0} and the
$p$th SB system \eqref{sb2.29}. In addition, assume
\eqref{sb3.1}--\eqref{sb3.2a} and let
$P\in\calK_{p}\setminus\{\Pinfp,\Pinfm,\Pzp,\Pzm\}$ and
$(n,n_{0})\in \bbZ^2$. Then for each $n\in \bbZ$, $\calD_{\humu(n)}$
and $\calD_{\hunu(n)}$ are nonspecial. Moreover,\footnote{To
avoid multi-valued expressions in formulas such as
\eqref{sb3.59}--\eqref{sb3.61b}, etc., we always agree to choose the same
path of integration connecting $Q_0$ and $P$.}
\begin{align}
\phi(P,n)&= C(n) \frac{\theta(\uz(P,\hunu(n)))}{\theta(\uz(P,\humu(n)))}
\exp\bigg(\int_{Q_0}^P \omega^{(3)}_{\Pzm, \Pinfm} \bigg), \lb{sb3.59} \\
\psi_{1}(P,n,n_{0})&=C(n,n_{0})\frac{\theta(\uz(P,\humu(n)))}
{\theta(\uz(P,\humu(n_{0})))} \exp\bigg((n-n_{0})
\int_{Q_{0}}^{P} \omega^{(3)}_{\Pzm,\Pinfp} \bigg), \lb{sb3.61} \\
\psi_{2}(P,n,n_{0})& =C(n) C(n,n_{0})  \no \\
& \hspace*{-1cm}\times
\frac{\theta(\uz(P,\hunu(n)))}{\theta(\uz(P,\humu(n_0)))}
\exp\bigg(\int_{Q_{0}}^{P}
\omega^{(3)}_{\Pzm,\Pinfm}+(n-n_{0}) \int_{Q_{0}}^{P}
\omega^{(3)}_{\Pzm,\Pinfp} \bigg),  \lb{sb3.61b}
\end{align}
where
\begin{align}
& C(n)=(-1)^{n-n_0}\exp\big[(n-n_0)(\omega_0^{0,-}(\Pzm,\Pinfm)
-\omega_0^{\infty_+}(\Pzm,\Pinfm))\big] \no \\
& \qquad \quad \, \times
\f{1}{\alpha(n_0)}\exp\big[-\omega_0^{0,+}(\Pzm,\Pinfm)\big]
\f{\theta(\uz(\Pzp,\humu(n_0)))}{\theta(\uz(\Pzp,\hunu(n_0)))},
\lb{sb3.61B} \\
&C(n,n_{0})= \exp\big[-(n-n_{0})\omega_0^{\infty_+}(\Pzm,\Pinfp)\big]
\frac{\theta(\uz(\Pinfp,\humu(n_{0})))}
{\theta(\uz(\Pinfp,\humu(n)))}. \lb{sb3.62}
\end{align}
The Abel map linearizes the auxiliary divisors in the sense that
\begin{align}
       \ual_{Q_{0}}(\calD_{\humu(n)})
&=\ual_{Q_{0}}(\calD_{\humu(n_{0})})
+\ua_{\Pzm}(\Pinfp)(n-n_{0}), \lb{sb3.63} \\
       \ual_{Q_{0}}(\calD_{\hunu(n)})
&=\ual_{Q_{0}}(\calD_{\hunu(n_{0})})
+\ua_{\Pzm}(\Pinfp)(n-n_{0}). \lb{sb3.64}
\end{align}
Finally, $\alpha,\beta$ are of the form
\begin{align}
&\alpha(n)=
\alpha(n_0)(-1)^{n-n_0}\exp\big[-(n-n_0)(\omega_0^{0,-}(\Pzm,\Pinfm)
-\omega_0^{\infty_+}(\Pzm,\Pinfm))\big] \no \\ & \hspace*{1.1cm} \times
\f{\theta(\uz(\Pzp,\hunu(n_0)))\theta(\uz(\Pzp,\humu(n)))}
{\theta(\uz(\Pzp,\humu(n_0)))\theta(\uz(\Pzp,\hunu(n)))}, \lb{sb3.66} \\
&\beta(n)= \beta(n_0) (-1)^{n-n_0}
\exp\big[(n-n_0)(\omega_0^{0,-}(\Pzm,\Pinfm)
-\omega_0^{\infty_+}(\Pzm,\Pinfm))\big] \no \\
& \hspace*{1.1cm} \times
\f{\theta(\uz(\Pinfp,\humu(n_0)))\theta(\uz(\Pinfp,\hunu(n)))}
{\theta(\uz(\Pinfp,\hunu(n_0)))\theta(\uz(\Pinfp,\humu(n)))}, \lb{sb3.68}
\\
& \alpha(n)\beta(n)=\exp\big[\omega_0^{\infty_+}(\Pzm,\Pinfm)
-\omega_0^{0,+}(\Pzm,\Pinfm)\big] \no \\
& \hspace*{1.8cm} \times
\f{\theta(\uz(\Pzp,\humu(n)))\theta(\uz(\Pinfp,\hunu(n)))}
{\theta(\uz(\Pzp,\hunu(n)))\theta(\uz(\Pinfp,\humu(n)))}. \lb{sb3.67}
\end{align}
\end{theorem}
\begin{proof}
While equation \eqref{sb3.63} is clear from \eqref{sb3.48}, equation  
\eqref{sb3.64} follows by combining \eqref{sb3.8} and \eqref{sb3.49}. By
Lemma \ref{sbl3.4}, $\calD_{\humu}$ and
$\calD_{\hunu}$ are nonspecial. By \eqref{sb3.8} and Theorem~\ref{taa17a},
$\phi(P,n)\exp\big(-\int_{Q_0}^P \omega^{(3)}_{\Pzm, \Pinfm} \big)$ must
be of the type  
\begin{equation}
              \phi(P,n)\exp\bigg(-\int_{Q_{0}}^{P}
              \omega^{(3)}_{\Pzm,\Pinfm} \bigg)= C(n)
\frac{\theta(\uz(P,\hunu(n)))}{\theta(\uz(P,\humu(n)))} \lb{sb3.69}
\end{equation}
for some constant $C(n)$. A comparison of \eqref{sb3.69} and the
asymptotic relations \eqref{sb3.46} and \eqref{sb3.47} then yields the
following expressions for $\alpha$ and $\beta$:
\begin{align}
(\alpha^+)^{-1}&= C^+ e^{\omega_0^{0,+}(\Pzm,\Pinfm)}
\f{\theta(\uz(\Pzp,\hunu^+))}{\theta(\uz(\Pzp,\humu^+))} \no \\
&=C^+ e^{\omega_0^{0,+}(\Pzm,\Pinfm)}
\f{\theta(\uz(\Pinfm,\hunu))}{\theta(\uz(\Pinfm,\humu))} \no \\
&=-C e^{\omega_0^{\infty_-}(\Pzm,\Pinfm)}
\f{\theta(\uz(\Pinfm,\hunu))}{\theta(\uz(\Pinfm,\humu))}. \lb{sb3.69a}
\end{align}
Similarly one obtains
\begin{align}
\beta^+&= C^+ e^{\omega_0^{\infty_+}(\Pzm,\Pinfm)}
\f{\theta(\uz(\Pinfp,\hunu^+))}{\theta(\uz(\Pinfp,\humu^+))} \no \\
&=C^+ e^{\omega_0^{\infty_+}(\Pzm,\Pinfm)}
\f{\theta(\uz(\Pzm,\hunu))}{\theta(\uz(\Pzm,\humu))} \no \\
&=-C e^{\omega_0^{0,-}(\Pzm,\Pinfm)}
\f{\theta(\uz(\Pzm,\hunu))}{\theta(\uz(\Pzm,\humu))}.  \lb{sb3.69d}
\end{align}
Here we used
\begin{equation}
\ual_{Q_0}(\calD_{\humu^+})=\ual_{Q_0}(\calD_{\humu})
+\ua_{\Pzm}(\Pinfp), \quad 
\ual_{Q_0}(\calD_{\hunu^+})=\ual_{Q_0}(\calD_{\hunu})
+\ua_{\Pzm}(\Pinfp),  \lb{sb3.69c}
\end{equation}
\eqref{sb3.58}, and relations of the type
\begin{align}
\uz(\Pinfp,\humu^+)&=\uz(\Pinfm,\hunu)=\uz(\Pzm,\humu)=\uz(\Pzp,\hunu^+),
\lb{sb3.72} \\
\uz(\Pinfp,\hunu^+)&=\uz(\Pzm,\hunu), \quad 
\uz(\Pzp,\humu^+)=\uz(\Pinfm,\humu). \lb{sb3.74} 
\end{align}
Thus, one concludes
\begin{equation}
C(n+1)=-\exp\big[\omega_0^{0,-}(\Pzm,\Pinfm)
-\omega_0^{\infty_+}(\Pzm,\Pinfm)\big]C(n), \quad n\in\bbZ  \lb{sb3.69e}
\end{equation}
and
\begin{equation}
C(n+1)=-\exp\big[\omega_0^{\infty_-}(\Pzm,\Pinfm)
-\omega_0^{0,+}(\Pzm,\Pinfm)\big]C(n), \quad n\in\bbZ,  \lb{sb3.69f}
\end{equation}
which is consistent with \eqref{sb3.57g}.
The first-order difference equation \eqref{sb3.69e} then implies
\begin{align}
C(n)&=(-1)^{(n-n_0)}\exp\big[(n-n_0)
(\omega_0^{0,-}(\Pzm,\Pinfm)
-\omega_0^{\infty_+}(\Pzm,\Pinfm))\big]C(n_0), \no \\
& \hspace*{8.5cm}  n,n_0\in\bbZ. \lb{sb3.69g}
\end{align}
Thus one infers \eqref{sb3.66} and \eqref{sb3.68}.
Moreover, \eqref{sb3.69g} and taking $n=n_0$ in the first line in
\eqref{sb3.69a} yield \eqref{sb3.61B}. Dividing the first line in
\eqref{sb3.69d} by the first line in \eqref{sb3.69a} then proves
\eqref{sb3.67}.

By \eqref{sb3.48} and Theorem \ref{taa17a}, $\psi_{1}(P,n,n_{0})$ must
be of the type \eqref{sb3.61}. A comparison of \eqref{sb3.11},
\eqref{sb3.46}, and \eqref{sb3.61} as $P\to\Pinfp$ ($\zeta=1/z$) then yields
\begin{align}
\psi_1(P,n,n_0)&\underset{\zeta\to 0}{=}\zeta^{n_0-n}(1+\Oh(\zeta)) 
\lb{sb3.74B}
\intertext{and}
\psi_1(P,n,n_0)&\underset{\zeta\to 0}{=}C(n,n_0)
\frac{\theta(\uz(\Pinfp,\humu(n)))}{\uz(\Pinfp,\humu(n_0)))} \no \\
& \quad \;\;\; \times
\exp\big[(n-n_0)\omega_0^{\infty_+}(\Pzm,\Pinfp)\big]
\zeta^{n_0-n}(1+\Oh(\zeta)) \lb{sb3.71}
\end{align}
proving \eqref{sb3.62}. Equation \eqref{sb3.61b} is clear from
\eqref{sb3.17}, \eqref{sb3.59}, and \eqref{sb3.61}.
\end{proof}

\begin{remark} \lb{sbremark3.7}
$(i)$ By \eqref{sb3.63}, \eqref{sb3.64}, the arguments of all theta
functions in \eqref{sb3.59}--\eqref{sb3.61b} \eqref{sb3.62}, and
\eqref{sb3.66}--\eqref{sb3.67} are linear with respect to $n$. \\
$(ii)$ Using relations of the type \eqref{sb3.72}, \eqref{sb3.74} and 
\begin{equation}
\ual_{Q_{0}}(\calD_{\hunu})
=\ual_{Q_{0}}(\calD_{\humu})+\ua_{\Pzm}(\Pinfm), \lb{sb3.75A}
\end{equation}
one can rewrite formulas \eqref{sb3.59}--\eqref{sb3.67} in terms of $\humu$
(or $\hunu$) only. \\
$(iii)$ For simplicity we assumed $\alpha(n)\neq 0$, $\beta(n)\neq 0$,
$n\in\bbZ$, in Theorem \ref{sbt3.6}. Since by \eqref{sb3.46} and
\eqref{sb3.47} no $\hmu_j$ and no $\hnu_\ell$ hits $\Pzp$ or $\Pinfp$, the
expressions \eqref{sb3.66} and \eqref{sb3.68} for $\alpha$ and $\beta$ are
consistent with this assumption. \\
$(iv)$ Generally, $\alpha$ and $\beta$ will not be quasi-periodic with
respect to $n\in\bbZ$. Only under certain restrictions on the
distribution of $\{E_m\}_{m=0,\dots,2p+1}$, such as the (de)focusing cases
discussed in Corollary \ref{sbc3.9} next, one can expect to uniformly bound
the exponential terms in \eqref{sb3.66} and \eqref{sb3.68} and prove
quasi-periodicity of $\alpha$ and $\beta$.  
\end{remark}

The special defocusing and focusing cases are briefly considered next.

\begin{corollary} \lb{sbc3.9}
Suppose that $\alpha, \beta \subset\bbC$ satisfy \eqref{sb3.0} and the
$p$th SB system \eqref{sb2.29} and assume
\eqref{sb3.1}--\eqref{sb3.2a}. Moreover, assume either
the defocusing case, where $\beta(n)=\ol{\alpha(n)}$, or the focusing case,
where $\beta(n)=-\ol{\alpha(n)}$, $n\in\bbZ$. In either case, $\alpha$ is
quasi-periodic with respect to $n\in\bbZ$.
\end{corollary}
\begin{proof}
We start by noting that the ratio of theta functions in \eqref{sb3.66} and
\eqref{sb3.68} is bounded as $n$ varies in $\bbZ$ since by \eqref{sb3.8}
(see also \eqref{sb3.46} and \eqref{sb3.47}) $\Pzp$ is never hit
by any $\hnu_\ell(n)$ and $\Pinfp$ is never hit by any $\hmu_j(n)$. Thus,
$\alpha$ (and of course $\beta$) is quasi-periodic if and only if the
exponential term in \eqref{sb3.66} is bounded (i.e., unimodular). Assume
the defocusing case $\beta=\ol\alpha$. Then, writing
\begin{equation}
\alpha(n)=b(n)e^{nc}, \quad \beta(n)=\tilde b(n)e^{-nc}, \; n\in\bbZ,
\quad b, \tilde b\in \ell^\infty(\bbZ)  \lb{3.75B}
\end{equation}
(cf.\ \eqref{sb3.66}, \eqref{sb3.68}), $\beta=\ol\alpha$ implies
\begin{equation}
\beta(n)=\tilde b(n) e^{-n\Re(c)-in\Im(c)}=\ol\alpha(n)=
\ol{b(n)} e^{n\Re(c)-in\Im(c)}  \lb{sb3.75C}
\end{equation}
and hence $\Re(c)=0$. The analogous argument applies in the focusing case.
\end{proof}

\begin{remark} \lb{sbr3.10}
$(i)$ The additional (de)focusing assumption $\beta=\pm\ol\alpha$ in
Corollary \ref{sbc3.9}, implies strong restrictions on the possible location
of the branch points $(E_m,0)$, $m=0,\dots,2p+1$. In particular, in analogy
to the Ablowitz--Ladik model discussed in
\cite{MillerErcolaniKricheverLevermore:1995}, one expects all $(E_m,0)$ to
occur in pairs which are reflection symmetric with respect to the unit
circle $\bbT$ in $\bbC$. In the defocusing case, $\beta=\ol\alpha$ with 
$|\alpha(n)|<1$, $n\in\bbZ$, all branch points are seen to lie on $\bbT$ as
discussed in \cite{GeronimoJohnson:1998} and \cite{Simon:2004}. For
$|\alpha|>1$ one expects them to bifurcate off the unit circle $\bbT$. \\
$(ii)$ In analogy to the defocusing case of the nonlinear
Schr\"odinger equation (cf.\ \cite[Ch.\ 3]{GesztesyHolden:2003}), the
isospectral manifold of algebro-geometric  solutions of \eqref{sb1.1} can
be identified with a $(p+1)$-dimensional real torus $\bbT^{p+1}$ as
discussed in detail in \cite[Ch.\ 11]{Simon:2004}. This isospectral torus is
of dimension $p+1$ (rather than $p$, given the $p$ divisors $\hmu_j(n_0)$,
$j=1,\dots,p$) due to the additional scaling invariance discussed in
\eqref{2.31}, \eqref{2.32} involving an arbitary constant multiple of
absolute value equal to one. \\
$(iii)$ By Remark \ref{sbr3.9}, no special divisors arise if
$\beta(n_1) =\pm\ol{\alpha(n_1)}=0$ for some $n_1\in\bbZ$
and hence Corollary \ref{sbc3.9} extends to this case as long as
$\beta(n_0)=\pm\ol{\alpha(n_0)} \neq 0$ in \eqref{sb3.68}. \\
$(iv)$ In the special defocusing case $\beta=\ol\alpha$, with 
$|\alpha(n)|<1$, $n\in\bbZ$, Corollary \ref{sbc3.9} recovers the original
result of Geronimo and Johnson \cite{GeronimoJohnson:1998} that $\alpha$ is
quasi-periodic without the use of Fay's generalized Jacobi variety, double
covers, etc.
\end{remark}

Finally, we briefly consider the case $p=0$ excluded in Theorem
\ref{sbt3.6}.

\begin{example}\lb{sbexample3.9}
Let $p=0$, $P=(z,y)\in\calK_0\backslash\{\Pzp,\Pzm,\Pinfp,\Pinfm\}$,
and $(n,n_0)\in\bbZ^2$. Then,
\begin{align}
&\calK_0 \colon \calF_0(z,y)=y^2-R_{2}(z)=y^2-(z-E_{0})(z-E_{1})=0,
\no \\
&  E_{0},E_{1}\in \bbC\backslash\{0\}, \; E_0\neq E_1, \quad g_1^2=E_0
E_1,
\quad g_1=y(\Pzp), \quad c_1=-(E_0+E_1)/2, \no \\
&\alpha(n)=\alpha(n_0)(-g_{1})^{n-n_0},  \quad 
\beta(n)=\beta(n_0)(-g_{1})^{n_0-n}, \no \\
& \sSB_{0}(\alpha,\beta)=
\begin{pmatrix} -2(\alpha^+ +g_1\alpha) \\
2(\beta^- +g_1\beta)\end{pmatrix}=0, \quad 
\alpha(n)\beta(n)=[1-(c_1/g_1)]/2, \no \\
&\phi(P)=\f{y+z-2\alpha^+\beta+c_1}{-2\alpha^+}
=\f{-2\beta z}{y-z+2\alpha^+\beta-c_1}. \no
\end{align}
\end{example}

One verifies that $E_0\neq E_1$ is equivalent to
$\alpha\beta\in\bbC\backslash\{0,1\}$. For a Borg-type theorem related to
this example in the special defocusing case $\beta=\ol\alpha$ with
$|\alpha(n)|<1$, $n\in\bbZ$, we refer to \cite{GesztesyZinchenko:2005}.

\appendix

\section{Hyperelliptic Curves and Their Theta Functions}\lb{A}
\renewcommand{\theequation}{A.\arabic{equation}}
\renewcommand{\thetheorem}{A.\arabic{theorem}}
\setcounter{theorem}{0}
\setcounter{equation}{0}

We give a brief summary of some of the fundamental properties
and notations needed from the
theory of hyperelliptic curves.  More details can be found in
some of the standard textbooks
\cite{FarkasKra:1992} and \cite{Mumford:1984}, as well as monographs
dedicated to integrable systems such as
\cite{BelokolosBobenkoEnolskiiItsMatveev:1994}, Ch.\ 2,
\cite{GesztesyHolden:2003}, App. A, B.

Fix $\N\in\bbN$. The hyperelliptic curve $\calK_\N$
of genus $\N$ used in Section~\ref{sbs3} is defined by
\begin{align}
&\calK_\N: \, \calF_\N(z,y)=y^2-R_{2\N+2}(z)=0, \quad
R_{2\N+2}(z)=\prod_{m=0}^{2\N+1}(z-E_m), \lb{b0} \\
& \{E_m\}_{m=0,\dots,2\N+1}\subset\bbC, \quad
E_m \neq E_{m'} \text{ for } m \neq m', \, m,m'=0,\dots,2\N+1.
\label{b1}
\end{align}
The curve \eqref{b1} is compactified by adding the
points $\Pinfp$ and $\Pinfm$,
$\Pinfp \neq \Pinfm$, at infinity.
One then introduces an appropriate set of
$\N+1$ nonintersecting cuts $\calC_j$ joining
$E_{m(j)}$ and $E_{m^\prime(j)}$. We denote
\begin{equation}
\calC=\bigcup_{j\in \{1,\dots,\N+1 \}}\calC_j,
\quad
\calC_j\cap\calC_k=\emptyset,
\quad j\neq k.\label{b2}
\end{equation}
Define the cut plane $\Pi=\bbC\backslash\calC$, 
and introduce the holomorphic function
\begin{equation}
R_{2\N+2}(\dott)^{1/2}\colon \Pi\to\bbC, \quad
z\mapsto \left(\prod_{m=0}^{2\N+1}(z-E_m) \right)^{1/2}\label{b4}
\end{equation}
on $\Pi$ with an appropriate choice of the square root
branch in \eqref{b4}. Define
\begin{equation}
\calM_{\N}=\{(z,\sigma R_{2\N+2}(z)^{1/2}) \mid
z\in\bbC,\; \sigma\in\{\pm 1\}
\}\cup\{\Pinfp,\Pinfm\} \label{b5}
\end{equation}
by extending $R_{2\N+2}(\dott)^{1/2}$ to $\calC$. The
hyperelliptic curve $\calK_\N$ is then the set
$\calM_{\N}$ with its natural complex structure obtained
upon gluing the two sheets of $\calM_{\N}$
crosswise along the cuts. The set of branch points
$\calB(\calK_\N)$ of $\calK_\N$ is given by
\begin{equation}
\calB(\calK_\N)=\{(E_m,0)\}_{m=0,\dots,2\N+1} \lb{5a}
\end{equation}
and finite points $P$ on $\calK_\N$ are denoted by
$P=(z,y)$, where $y(P)$ denotes the meromorphic function
on $\calK_\N$ satisfying $\calF_\N(z,y)=y^2-R_{2\N+2}(z)=0$.
Local coordinates near $P_0=(z_0,y_0)\in\calK_\N\backslash
(\calB(\calK_\N)\cup\{\Pinfp,\Pinfm\})$ are
given by $\zeta_{P_0}=z-z_0$, near $\Pinfpm$ by
$\zeta_{\Pinfpm}=1/z$, and near branch points
$(E_{m_0},0)\in\calB(\calK_\N)$ by
$\zeta_{(E_{m_0},0)}=(z-E_{m_0})^{1/2}$. The Riemann surface
$\calK_\N$ defined in this manner has topological genus $\N$.

One verifies that $dz/y$ is a holomorphic differential
on $\calK_\N$ with zeros of order $\N-1$ at $\Pinfpm$
and that
\begin{equation}
\eta_j=\frac{z^{j-1}dz}{y}, \quad j=1,\dots,\N
\lb{b24}
\end{equation}
form a basis for the space of holomorphic differentials
on $\calK_\N$.  Introducing the
invertible matrix $C$ in $\bbC^\N$,
\begin{align}
\begin{split}
C & =(C_{j,k})_{j,k=1,\dots,\N}, \quad C_{j,k}
= \int_{a_k} \eta_j, \\
\underline{c} (k) & = (c_1(k), \dots,
c_n(k)), \quad c_j (k) =
C_{j,k}^{-1}, \;\, j,k=1,\dots,g, \lb{A.7}
\end{split}
\end{align}
the corresponding basis of normalized holomorphic
differentials $\omega_j$, $j=1,\dots,\N$ on $\calK_\N$ is given by
\begin{equation}
\omega_j = \sum_{\ell=1}^\N c_j (\ell) \eta_\ell,
\quad \int_{a_k} \omega_j =
\delta_{j,k}, \quad j,k=1,\dots,\N. \lb{b26}
\end{equation}
Here $\{a_j,b_j\}_{j=1,\dots,\N}$ is a homology basis for
$\calK_\N$ with intersection matrix of the cycles satisfying
\begin{equation}
a_j \circ b_k=\delta_{j,k}, \; a_j \circ a_k=0,
\; b_j \circ b_k=0, \quad j,k=1,\dots,\N. \lb{c26}
\end{equation}

Associated with the homology basis
$\{a_j, b_j\}_{j=1,\dots,\N}$ we
also recall the canonical dissection of $\calK_\N$
along its cycles yielding
the simply connected interior $\hatt \calK_\N$ of the
fundamental polygon $\partial {\hatt \calK}_\N$ given by
$\partial  {\hatt \calK}_\N =a_1 b_1 a_1^{-1} b_1^{-1}
a_2 b_2 a_2^{-1} b_2^{-1} \cdots a_\N^{-1} b_\N^{-1}$.
Let $\calM (\calK_\N)$ and $\calM^1 (\calK_\N)$ denote the
set of meromorphic
functions (0-forms) and meromorphic
differentials (1-forms)
on $\calK_\N$. The residue of a meromorphic differential
$\nu\in \calM^1 (\calK_\N)$ at a
point $Q \in \calK_\N$ is defined by $\text{res}_{Q}(\nu)
=\frac{1}{2\pi i} \int_{\gamma_{Q}} \nu$, 
where $\gamma_{Q}$ is a counterclockwise oriented, smooth, simple, closed
contour encircling $Q$ but no other pole of
$\nu$.  Holomorphic
differentials are also called Abelian differentials
of the first kind. Abelian differentials of the
second kind, $\omega^{(2)} \in \calM^1 (\calK_\N)$, are characterized
by the property that all their residues vanish.  Any meromorphic
differential $\omega^{(3)}$ on $\calK_\N$ not of the first or second kind is
said to be of the third kind. A differential of the third kind
$\omega^{(3)} \in \calM^1 (\calK_\N)$ is usually normalized by the
vanishing of its $a$-periods, that is,
\begin{equation}
\int_{a_j} \omega^{(3)} =0, \quad  j=1,\dots, \N.
\lb{a36}
\end{equation}
A normal differential of the third kind $\omega_{P_1, P_2}^{(3)}$ associated
with two points $P_1$,
$P_2 \in \hatt \calK_\N$, $P_1 \neq P_2$ by definition
has simple poles at
$P_j$ with residues $(-1)^{j+1}$, $j=1,2$ and
vanishing $a$-periods.  If $\omega_{P,Q}^{(3)}$ is a
normal differential of the third kind associated
with $P$, $Q\in\hatt \calK_\N$, holomorphic on
$\calK_\N \backslash \{ P,Q\}$, then
\begin{equation}
\frac{1}{2\pi i} \int_{b_j} \omega_{P,Q}^{(3)} = \int_{Q}^P \omega_j,
\quad  j=1,\dots,\N,
\lb{a37}
\end{equation}
where the path from $Q$ to $P$ lies in
$\hatt \calK_\N$ (i.e.,
does not touch any of the cycles $a_j$, $b_j$). 

We shall always assume (without loss of generality)
that all poles of differentials of the second and third kind
on $\calK_\N$ lie on $\hatt \calK_\N$ (i.e.,
not on $\partial \hatt \calK_n$).

Define the matrix $\tau=(\tau_{j,\ell})_{j,\ell=1,\dots,\N}$ by
\begin{equation}
\tau_{j,\ell}=\int_{b_\ell}\omega_j, \quad j,\ell=1,
\dots,\N. \label{b8}
\end{equation}
Then $\Im(\tau)>0$ and $\tau_{j,\ell}=\tau_{\ell,j}$,
$j,\ell =1,\dots,\N$. Associated with $\tau$ one introduces the period
lattice
\begin{equation}
L_\N = \{ \ul z \in\bbC^\N \mid \ul z = \ul m + \ul n\tau,
\; \ul m, \ul n \in\bbZ^\N\}
\lb{a28}
\end{equation}
and the Riemann theta function associated with $\calK_\N$ and
the given homology basis $\{a_j,b_j\}_{j=1,\dots,\N}$,
\begin{equation}
\theta(\ul z)=\sum_{\ul n\in\bbZ^\N}\exp\big(2\pi
i(\ul n,\ul z)+\pi i(\ul n, \ul n\tau)\big),
\quad \ul z\in\bbC^\N, \label{b9}
\end{equation}
where $(\ul u, \ul v)=\ol {\ul u}\, \ul v^\top=\sum_{j=1}^\N \ol{u_j}\, v_j$
denotes the scalar product in $\bbC^\N$. It has the fundamental properties
\begin{align}
& \theta(z_1, \ldots, z_{j-1}, -z_j, z_{j+1},
\ldots, z_\N) =\theta
(\ul z), \lb{a27}\\
& \theta (\ul z +\ul m + \ul n\tau)
=\exp \big(-2 \pi i (\ul n,\ul z) -\pi i (\ul n, \ul n\tau) \big) \theta
(\ul z), \quad \ul m, \ul n \in\bbZ^\N.
\lb{aa51}
\end{align}

Next, fix a base point $Q_0\in\calK_\N\backslash
\{\Pzpm,\Pinfpm\}$, denote by
$J(\calK_\N) = \bbC^\N/L_\N$ the Jacobi variety of $\calK_\N$,
and define the
Abel map $\underline{A}_{Q_0}$ by
\begin{equation}
\underline{A}_{Q_0} \colon \calK_\N \to J(\calK_\N), \quad
\underline{A}_{Q_0}(P)=
\bigg(\int_{Q_0}^P \omega_1,\dots,\int_{Q_0}^P \omega_\N \bigg)
\pmod{L_\N}, \quad P\in\calK_\N. \label{b10}
\end{equation}
Similarly, we introduce
\begin{equation}
\ul \alpha_{Q_0}  \colon
\Div(\calK_\N) \to J(\calK_\N),\quad
\calD \mapsto \ul \alpha_{Q_0} (\calD)
=\sum_{P \in \calK_\N} \calD (P) \ul A_{Q_0} (P),
\label{aa47}
\end{equation}
where $\Div(\calK_\N)$ denotes the set of
divisors on $\calK_\N$. Here $\calD \colon \calK_\N \to \bbZ$
is called a divisor on $\calK_\N$ if $\calD(P)\neq0$ for only
finitely many $P\in\calK_\N$. (In the main body of this paper
we will choose $Q_0$ to be one of the branch points, i.e.,
$Q_0\in\calB(\calK_\N)$, and for simplicity we will always choose
the same path of integration from $Q_0$ to $P$ in all Abelian
integrals.) 

In connection with divisors on $\calK_\N$ we shall employ the
following
(additive) notation,
\begin{align} \lb{A.17}
&\calD_{Q_0\ul Q}=\calD_{Q_0}+\calD_{\ul Q}, \quad \calD_{\ul
Q}=\calD_{Q_1}+\cdots +\calD_{Q_m}, \\
& {\ul Q}=\{Q_1, \dots ,Q_m\} \in \sym^m \calK_\N,
\quad Q_0\in\calK_\N, \; m\in\bbN, \no
\end{align}
where for any $Q\in\calK_\N$,
\begin{equation} \lb{A.18}
\calD_Q \colon  \calK_\N \to\bbN_0, \quad
P \mapsto  \calD_Q (P)=
\begin{cases} 1 & \text{for $P=Q$},\\
0 & \text{for $P\in \calK_\N\backslash \{Q\}$}, \end{cases}
\end{equation}
and $\sym^n \calK_\N$ denotes the $n$th symmetric product of
$\calK_\N$. In particular, $\sym^m \calK_\N$ can be
identified with
the set of nonnegative
divisors $0 \leq \calD \in \Div(\calK_\N)$ of degree $m$.

For $f\in \calM (\calK_\N) \backslash \{0\}$,
$\omega \in \calM^1 (\calK_\N) \backslash \{0\}$ the
divisors of $f$ and $\omega$ are denoted
by $(f)$ and
$(\omega)$, respectively.  Two
divisors $\calD$, $\calE\in \Div(\calK_\N)$ are
called equivalent, denoted by
$\calD \sim \calE$, if and only if $\calD -\calE
=(f)$ for some
$f\in\calM (\calK_\N) \backslash \{0\}$.  The divisor class
$[\calD]$ of $\calD$ is
then given by $[\calD]
=\{\calE \in \Div(\calK_\N)\mid\calE \sim \calD\}$.  We
recall that
\begin{equation}
\deg ((f))=0,\, \deg ((\omega)) =2(\N-1),\,
f\in\calM (\calK_\N) \backslash
\{0\},\,  \omega\in \calM^1 (\calK_\N) \backslash \{0\},
\lb{a38}
\end{equation}
where the degree $\deg (\calD)$ of $\calD$ is given
by $\deg (\calD)=\sum_{P\in \calK_\N} \calD (P)$. 
$(f)$ is called a principal divisor.

Introducing the complex linear spaces
\begin{align}
\calL (\calD) & =\{f\in \calM (\calK_\N)\mid f=0
       \text{ or } (f) \geq \calD\}, \quad
r(\calD) =\dim \calL (\calD),
\lb{a39}\\
\calL^1 (\calD) & =
       \{ \omega\in \calM^1 (\calK_\N)\mid \omega=0
       \text{ or } (\omega) \geq
\calD\}, \quad i(\calD) =\dim \calL^1 (\calD) \lb{a40}
\end{align}
with $i(\calD)$ the index of speciality of $\calD$, one
infers that $\deg(\calD)$, $r(\calD)$, and $i(\calD)$ only depend on
the divisor class $[\calD]$ of $\calD$.  Moreover, we recall the
following fundamental facts.

\begin{theorem} \lb{thm1}
Let $\calD \in \Div(\calK_\N)$,
$\omega \in \calM^1 (\calK_\N) \setminus \{0\}$. Then
\begin{equation}
       i(\calD) =r(\calD-(\omega)), \quad \N\in\bbN_0.
\lb{a41}
\end{equation}
The Riemann--Roch theorem reads
\begin{equation}
r(-\calD) =\deg (\calD) + i (\calD) -\N+1,
\quad \N\in\bbN_0.
\lb{a42}
\end{equation}
By Abel's theorem, $\calD\in \Div(\calK_\N)$,
$\N\in\bbN$, is principal
if and only if
\begin{equation}
\deg (\calD) =0 \text{ and } \ul \alpha_{Q_0} (\calD)
=\ul{0}.
\lb{a43}
\end{equation}
Finally, assume
$\N\in\bbN$. Then $\ul \alpha_{Q_0}
: \Div(\calK_\N) \to J(\calK_\N)$ is surjective
$($Jacobi's inversion theorem$)$.
\end{theorem}

\begin{theorem} \lb{thm3}
Let $\calD_{\ul Q} \in \sym^\N \calK_\N$,
$\ul Q=\{Q_1, \ldots, Q_\N\}$.  Then $1 \leq i (\calD_{\ul Q} ) =s \leq
\N/2$ if and only if there are $s$ pairs of the type
$(P, P^*)\in \{Q_1,\ldots, Q_\N\}$ $($this includes, of course, branch
points for which $P=P^*$$)$.
\end{theorem}

Denote by $\ul \Xi_{Q_0}=(\Xi_{Q_{0,1}}, \dots,
\Xi_{Q_{0,\N}})$ the vector of Riemann constants,
\begin{equation}
\Xi_{Q_{0,j}}=\frac12(1+\tau_{j,j})-
\sum_{\substack{\ell=1 \\ \ell\neq j}}^\N\int_{a_\ell}
\omega_\ell(P)\int_{Q_0}^P\omega_j,
\quad j=1,\dots,\N. \lb{aa55}
\end{equation}

\begin{theorem} \lb{taa17a}
Let $\ul Q =\{Q_1,\dots,Q_\N\}\in \sym^\N \calK_\N$ and
assume $\calD_{\ul Q}$ to be nonspecial, that is,
$i(\calD_{\ul Q})=0$. Then
\begin{equation}
\theta(\ul {\Xi}_{Q_0} -\ul {A}_{Q_0}(P) + \alpha_{Q_0}
(\calD_{\ul Q}))=0 \text{ if and only if }
P\in\{Q_1,\dots,Q_\N\}. \lb{aa55a}
\end{equation}
\end{theorem}
 
\medskip

{\bf Acknowledgments.}
We thank Russell Johnson and Barry Simon for discussions. F.G. gratefully
acknowledges the extraordinary hospitality of the Department of
Mathematical Sciences of the Norwegian University of Science and
Technology, Trondheim, during extended stays in the summers of 2001--2004,
where parts of this paper were written.


\end{document}